\documentclass[]{article}
\usepackage{graphicx,subfigure}
\usepackage{amsfonts,amsbsy,amscd,amsgen,amsthm,dsfont,amsmath,amssymb,mathrsfs}
\usepackage[colorlinks=true,allcolors=blue]{hyperref}
\usepackage{bm,multirow,esvect}
\usepackage{tabularx}
\usepackage{enumitem}
\usepackage{xcolor}
\usepackage{lipsum}
\usepackage{authblk}

\newcommand{\bbeta}{\pmb \beta}
\newcommand{\btheta}{\pmb\theta}

\newcommand{\id}{\mathrm d}
\newcommand{\bxi}{\pmb \xi}
\newcommand{\vc}{\mathbf}

\newcommand{\etab}{\pmb{\eta}}

\renewcommand{\i}{\hat{i}}
\newcommand{\pard}[2]{\frac{\partial #1}{\partial #2}}

\newcommand{\B}{\mathcal{B}}
\renewcommand{\H}{\mathcal{H}}
\newcommand{\g}{\mathrm{g}}
\newcommand{\M}{\mathcal{M}}

\DeclareMathAlphabet\mathbfcal{OMS}{cmsy}{b}{n}

\newcommand{\spn}{\mbox{span}}

\theoremstyle{definition}
\newtheorem{thm}{Theorem}
\newtheorem{rem}{Remark}

\newtheorem{defn}{Definition}
\newtheorem{ass}{Assumption}

\begin{document}
\title{A review of shape-morphing solutions and evolutional neural networks for spatiotemporal dynamics} 
\author{Mohammad Farazmand\thanks{Corresponding author: farazmand@ncsu.edu}}
\affil[]{Department of Mathematics, North Carolina State University, 2311 Stinson Drive, Raleigh, NC 27695-8205, USA}
\date{}
\maketitle

\begin{abstract}
Shape-morphing solutions (SMS) refer to a class of approximate solutions of partial differential equations (PDEs) with the distinguishing feature that they depend nonlinearly on a set of time-dependent parameters. They generalize Galerkin truncations by allowing the basis (or trial functions) to evolve in time in order to adapt to the solution of the PDE. As such, SMS are particularly suitable for reduced-order modeling as well as high fidelity simulation of multiscale systems which exhibit localized time-dependent features, such as vortices, dispersive wave packets, and shocks. Furthermore, being mesh-free, SMS is scalable for solving PDEs in higher spatial dimensions. As a special case, SMS allows the approximation of the PDE's solution by a neural network whose weights and biases depend on time. Such neural networks are known as evolutional neural networks or neural Galerkin schemes. The evolution of SMS parameters is dictated by the SMS equation, a set of ordinary differential equations derived from the Dirac--Frenkel variational principle.  Over the past five years, contributions to the theory and computation of SMS have been growing rapidly. Here, we survey these developments, showcase some applications of SMS, and highlight important open problems for future research. At the same time, this review is structured to serve as a tutorial for applied mathematicians, physicist, and engineers who wish to enter this field. 
\end{abstract}

\tableofcontents

\section{Introduction}\label{sec:intro}
Consider a general evolutionary partial differential equation (PDE),
\begin{equation}
	\partial_t u = F(u), \quad u(\vc x,0)=u_0(\vc x),
	\label{eq:gen_PDE}
\end{equation}
where $F$ is a potentially nonlinear differential operator. We assume that this initial value problem is well-posed and has a unique solution (in an appropriate sense) for all times $t\geq 0$. A conventional Galerkin-type approximation of the solution $u(\vc x,t)$ takes the form 
\begin{equation}
u(\vc x,t) \simeq \sum_{i=1}^r \alpha_i(t) \phi_i(\vc x),
\label{eq:Galerkin}
\end{equation}
where $\alpha_i(t)$ are the time-dependent parameters whose evolution needs to be determined. The choice of the basis (or trial) functions $\phi_i(\vc x)$, which are independent of time, results in different numerical methods. For example, choosing Fourier modes as the basis leads to the Fourier (pseudo-) spectral method and piecewise linear (hat) functions lead to the finite element method.
\begin{figure}[t!]
	\centering
	\includegraphics[width=.95\textwidth]{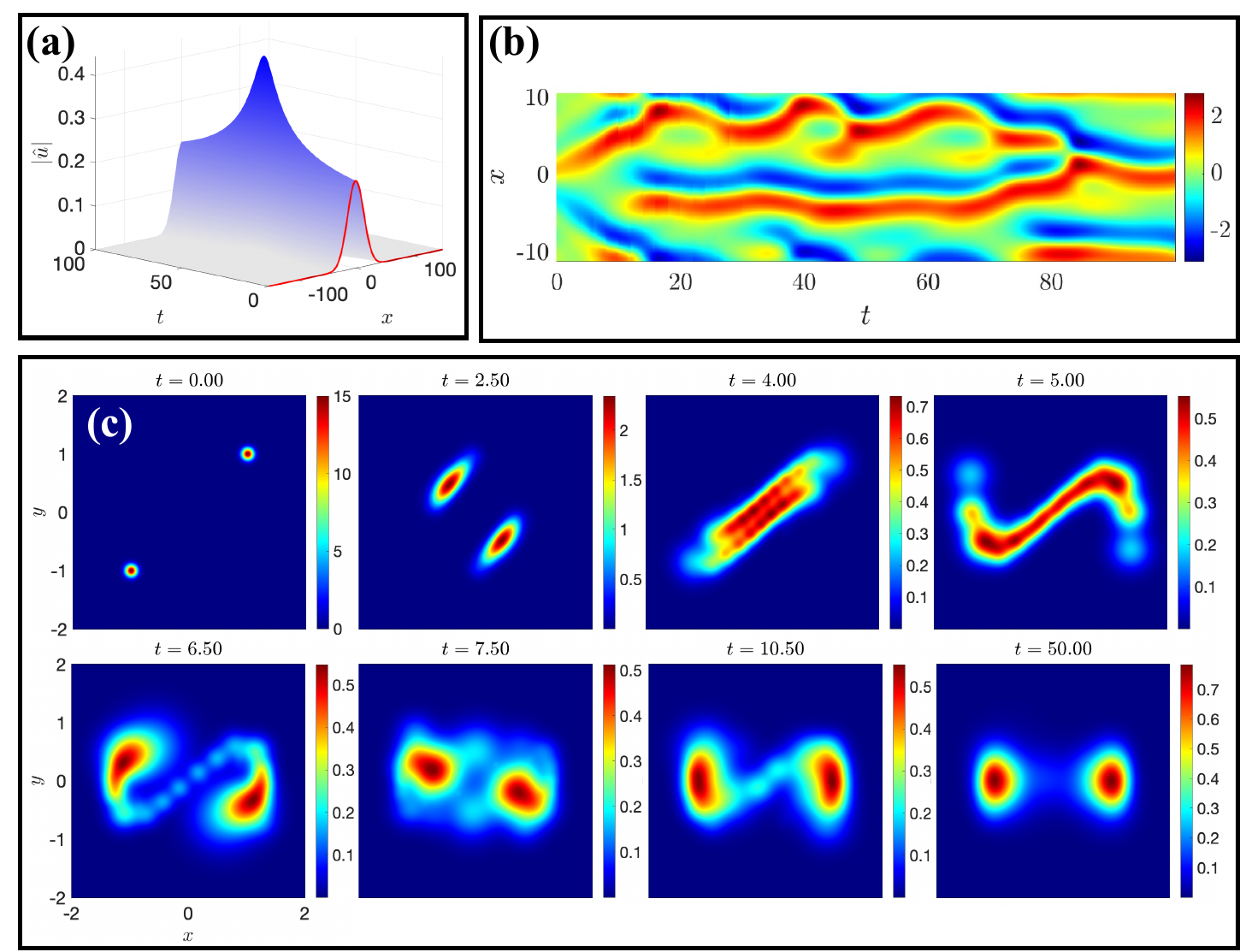}
	\caption{Three examples of shape-morphing solutions solving (a) the nonlinear Schr\"odinger equation, (b) the Kuramoto--Sivashinsky equation, (c) the Fokker--Planck equation. These results are discussed in detail in Section~\ref{sec:num_ex}.}
	\label{fig:num_ex}
\end{figure}

However, the static nature of the basis $\phi_i$ may lead to computational inefficiencies when the solution comprises time-varying localized structures. Examples include vortices in fluid dynamics,  shock waves, jet streams in the atmosphere, and dispersive wave packets. For such PDEs, an unnecessarily large number of static basis functions may be needed to reproduce a localized, yet time-varying, feature. 

Shape-morphing solutions (SMS) address this issue by allowing the basis functions to also evolve over time, hence moving and deforming in order to adapt to the solution of the PDE (see figure \ref{fig:num_ex}, for examples). As a special case, consider the SMS,
\begin{equation}
\hat u(\vc x,\btheta(t)):= \sum_{i=1}^r \alpha_i(t) \phi_i(\vc x,\bbeta_i(t)),
\label{eq:sms_Galerkin}
\end{equation}
where in addition to the amplitudes $\alpha_i(t)\in\mathbb R$, each mode depends on the \emph{shape parameters} $\bbeta_i(t)\in\mathbb R^{K-1}$ which encode the time-dependence of the basis functions. Together, the amplitudes and shape parameters form the set of all parameters $\btheta(t) = [\alpha_1,\bbeta_1^\top,\cdots, \alpha_r,\bbeta_r^\top]\in\mathbb R^n$ where $n=rK$.
In Section~\ref{sec:sms_ex}, we provide several concrete examples of SMS in addition to the working example described in Section~\ref{sec:working}.

In applying SMS to a given PDE, several theoretical and computational questions immediately arise. 
For example, (i) what is the optimal shape of the SMS $\hat u(\vc x,\btheta(t))$ and its dependence on the parameters $\btheta$? 
(ii) How does the approximation error depend on the shape of the SMS and the number of parameters $n$?
(iii) How should the parameters $\btheta(t)$ evolve so that an SMS $\hat u(\vc x,\btheta(t))$ best approximates the solution of the PDE? 
(iv) What computational methods should be used for efficient implementation of SMS? As we review here, some of these questions have been addressed (at least partially) and some remain open problems.

The idea of SMS was introduced by Anderson and Farazmand~\cite{Anderson2022} in the context of reduced-order modeling (ROM) and independently by Du and Zaki~\cite{Zaki2021} in the context of neural networks. When SMS are used for reduced-order modeling, they are also referred to as Reduced-Order Nonlinear Solutions (RONS). 
When used as a neural network, they are called an Evolutional Deep Neural Network (EDNN).

Ever since their introduction, there has been a fast growing body of work studying the theoretical and computational aspects of SMS. This review paper seeks to survey these developments and highlight potential areas for future research.

\subsection{Historical background}\label{sec:history}
Studies that are directly related to SMS are reviewed throughout this paper. In this section, we take the opportunity to give a historical overview of more conventional methods that seek to address the same underlying problems. In particular, we review adaptive mesh refinement, multigrid methods, and mesh-free solvers. Readers who are not interested in this historical background may safely skip this section.

\emph{Adaptive methods}:
Adaptive methods use an evolving grid to adequately resolve all scales. A smaller grid size is used to resolve localized small-scale structures, whereas the large-scale surroundings are resolved using a larger grid size. Since small scale structures deform and move, the grid needs to be updated during time stepping. Adaptive methods are ubiquitously used in various scientific fields where multiscale PDEs arise, such as in ocean modeling~\cite{Blayo1999,Cornford2013,leveque_2011}, meteorology~\cite{Ferguson2016,Krishnamurti1995,Toth2012}, astrophysics~\cite{Klein1999}, and fluid-structure interactions~\cite{Bayyuk1993,Berger1989,Berger1985,vanella2010}.

A naive implementation of adaptive methods can be prohibitive since mesh generation at every time step is computationally expensive.
Furthermore, special attention needs to be given to the time stepping scheme itself, especially for hyperbolic PDEs, so that it is stable across all grid sizes. These considerations have been addressed in a class of algorithms referred to as Adaptive Mesh Refinement (AMR) methods which systematically refine the mesh in the regions where higher resolution is needed~\cite{Berger1984}. 
Adaptive grids can be used in conjunction with familiar discretization methods such as finite difference~\cite{Berger1989,Berger1998,Berger1984}, finite element~\cite{Babuska1978,Bank1981}, and finite volume~\cite{Cornford2013,Jasak2000,Zhang2012} schemes. More recently, even higher-order spectral/hp element methods have been used on adaptive grids~\cite{tanarro2020}.

Another class of adaptive methods, primarily used in fluid dynamics,  is the so-called 
arbitrary Lagrangian-Eulerian (ALE) method~\cite{HIRT1974,Hirt1997,Hirt1970,Trulio1966}. The main idea behind the ALE method
is to move the grid with the fluid with the primary purpose of resolving localized small-scale structures and shock waves. This requires tracking the Lagrangian (or material) deformation of the fluid 
on an evolving grid in addition to its evolution on a fixed Eulerian grid.
Although ALE is often used with finite difference or finite volume discretizations, it has also been adopted for use with spectral/hp element methods~\cite{Helenbroook2001,Helenbrook2002,karniadakis2013}

We point out that the idea of changing the basis function is already implicit in adaptive methods. The evolution of the mesh can be recast as changing the basis over which the PDE solution is being approximated. Shape-morphing solutions make this notion explicit by optimally and automatically changing the trial functions as opposed to the grid.

\emph{Multiresolution and multigrid methods}:
Wavelet methods were originally developed for signal processing and data compression~\cite{Daubechies1988,gabor1946,haar1910,Kronland1987,meyer1986,zweig1976}. 
Unlike the Fourier basis which is local only in the frequency domain, wavelets can be local in both space and frequency domains.
It was quickly realized that wavelets can be used for analysis of multiscale phenomena beyond signal processing. For instance, wavelets were used to analyze coherent structures in two- and three-dimensional turbulence~\cite{kevlahan07,farge_3d,farge}. Wavelet-based methods have also been developed to numerically solve multiscale PDEs, both with Galerkin truncation~\cite{Bertoluzza1995,Beylkin1991,BEYLKLN1997,DAHMEN2001,dahmen1993} and with collocation points~\cite{Bertoluzza1997,kevlahan2005,schneider1997}.
Multiresolution analysis (MRA) provides a systematic method for constructing wavelets which from a basis for $L^2$, and more general Hilbert or Banach spaces~\cite{Daubechies1988,Walnut1989,Mallat1989,meyer1986,Walnut2001}. Derivation based on MRA yields a hierarchy of subspaces which correspond to ever smaller length scales that are localized in space. Note the resemblance to multigrid methods~\cite{brandt1973,briggs1988,Fedorenko1962,McCormick1989,Owhadi2017} which encode various scales in a hierarchy of grids
as opposed to basis functions.

It is important to note that wavelet methods still use a time-independent basis (or set of modes) for the function space. 
The number of wavelets grows exponentially with the resolution level ($2^j$ for level $j$). Consequently, wavelet methods are quite memory intensive. Adaptive wavelet methods have been introduced to manage these memory issues~\cite{Deiterding2020,kevlahan2005,Mehra2008b,Mehra2008a,Vasilyev2000}. 
In contrast, SMS evolves the basis over time to adapt to the solution of the PDE. As a result, far fewer modes are required alleviating the memory issues associated with wavelets.

\emph{Mesh-free solvers}: 
To avoid the computationally expensive steps of generating and updating a computational grid, a variety of mesh-free methods have been developed. Examples include the smoothed particle hydrodynamics (SPH)~\cite{Monaghan1977,Lucy1977}, radial basis functions~\cite{Kansa1990b}, diffuse element method \cite{Nayroles1992}, reproducing kernel particle method~\cite{Liu1995}, and meshless local Petrov--Galerkin method \cite{Atluri1998}. Perhaps the most prominent among these methods is SPH which estimates a field quantity as a convolution with a smooth kernel. These kernels can be interpreted as a weighted average of the contribution of each fluid particle to the field quantity of interest. The PDE dictates the distribution of these particles and their associated kernels. We observe that certain shape-morphing solutions can be viewed as a mesh-free method; for instance, see examples~\eqref{eq:sms_gauss_gen} and~\eqref{eq:sms_nn} in Section~\ref{sec:sms_ex} below.

\emph{Dynamical low-rank approximation}: 
Digressing somewhat from the subject of PDEs, we note that the idea of adaptive and time-dependent subspaces appears in the literature on matrix differential equations, $\dot A = F(A)$, where $A(t)$ is a time-dependent matrix and $F(A)$ is a matrix-valued map. Koch and Lubich~\cite{Koch2007} proposed the idea of Dynamical Low-Rank Approximation (DLRA) to evolve the matrix differential equations with low-rank updates to the matrix $A(t)$; also see Ref.~\cite{Wiesel1993}. When the time-derivative $\dot A$ is a low-rank matrix, this produces accurate solutions with a relatively low computational cost.
DLRA accomplishes this by adaptively changing the low-dimensional subspace to which the updates belong. A similar idea appears in the Optimally Time-Dependent (OTD) basis of Babaee and Sapsis~\cite{otd} who proposed a low-rank approximation of instabilities for linearized PDEs. Both DLRA and OTD use
a Dirac--Frenkel variational principle to derive the evolution of their low-dimensional subspaces. Although the scope and aims of DLRA/OTD are clearly different from SMS, the underlying variational principles are similar. We further comment on this in Remark~\ref{rem:DF_vari}.

\subsection{Outline}
The remainder of this review is organized as follows: 
\begin{itemize}[wide,labelwidth=!,labelindent=0pt]
	\item Section~\ref{sec:working} introduces a simple working example on which the SMS theory is demonstrated. 
	\item Section~\ref{sec:theory} reviews the main theoretical results.
	\item Section~\ref{sec:geom} described a geometric interpretation of SMS and its implications.
	\item Section~\ref{sec:EDNN} elucidates the connections between SMS and evolutional neural networks.
	\item Section~\ref{sec:num_ex} presents four numerical examples, showcasing the applications of SMS.
	\item Section~\ref{sec:comp} discusses the computational aspects of SMS.
	\item Section~\ref{sec:concl} presents concluding remarks and future directions.
\end{itemize}

\section{Working example}\label{sec:working}
As we review in this paper, the SMS framework is applicable to general, nonlinear and challenging PDEs. 
However, it is helpful to consider a simple working example on which the theory and computational methods are illustrated. For such a working example, we choose the linear advection-diffusion equation,
\begin{equation}\label{eq:ADE}
	u_t + cu_x = \nu u_{xx}, \quad u(x,0)=\delta(x),
\end{equation}
where $c$ is the advection speed, $\nu>0$ is the diffusion coefficient, and $\delta$ denotes the Dirac delta measure. For this PDE, the differential operator $F$ is linear, given by $F(u) = - cu_x + \nu u_{xx}$.
The fundamental solution of~\eqref{eq:ADE} is given by 
\begin{equation}\label{eq:ADE_FS}
	u(x,t) = \frac{1}{\sqrt{4\pi \nu t}}\exp\left[- \frac{(x-ct)^2}{4\nu t}\right], 
\end{equation}
which is smooth for all $t>0$.

On the other hand, we consider an SMS, 
\begin{equation}\label{eq:ADE_sms}
	\hat u(x,\btheta(t))=\alpha(t)\exp \left[ - \left(\frac{x-\mu(t)}{\sigma(t)}\right)^2\right].
\end{equation}
comprising a single Gaussian with time-dependent amplitude $\alpha(t)$, mean $\mu(t)$, and lengthscale $\sigma(t)$. The shape parameters are given by $\bbeta(t) = [\mu(t),\sigma(t)]^\top$, resulting in the complete set of parameters $\btheta(t) = [\alpha(t),\mu(t),\sigma(t)]^\top$. The parameter $\mu(t)$ allows the Gaussian to move in space, whereas $\sigma(t)$ allows for it to change its lengthscale, or `breath'.

With the correct evolution of parameters $\btheta(t)\in\mathbb R^3$, the SMS~\eqref{eq:ADE_sms} coincides with the fundamental solution~\eqref{eq:ADE_FS}. 
We will see later, in Section~\ref{sec:sms_eq}, that the SMS equations indeed return the correct evolution of parameters; see Eq.~\eqref{eq:qdot_we}. Note that here we chose the SMS shape~\eqref{eq:ADE_sms} to resemble the fundamental solution. In general, an appropriate shape for the SMS depends on the PDE and requires some familiarity with typical solutions of it. We further comment on the optimal choice of the SMS in Section~\ref{sec:concl}.

Throughout this paper, we use this working example to illustrate theoretical and computational aspects of SMS. 
Carrying out the required calculations manually could be laborious. However, using a symbolic computing software such as Mathematica or Maple, the computations can be carried out within seconds.

\section{Shape-morphing solutions}\label{sec:theory}
The PDE~\eqref{eq:gen_PDE} is defined on a $d$-dimensional spatial domain $\vc x\in\mathcal D \subseteq \mathbb R^d$ with the appropriate boundary conditions. The solution to the PDE is therefore a function $u:\mathcal D \times \mathbb R^+\to\mathbb R^p$. In the working example~\eqref{eq:ADE}, we have $d=1$, $\mathcal D=\mathbb R$, and the boundary conditions only require that the solution decays to zero at infinity, $\lim_{|\vc x|\to\infty} u(\vc x,t)=0$.

In this review, we only consider scalar-valued PDEs ($p=1$), but generalization to higher dimensions is straightforward; see, e.g.,~\cite{Hilliard2024}. Furthermore, we only consider strong solutions of the PDE. Generalization to weak solutions is not so straightforward.

We emphasize that, although we consider scalar-valued PDEs, the spatial domain $\mathcal D$ is allowed to have a higher dimension $d$. In fact, the computational efficiency of SMS, compared to conventional methods, becomes more pronounced as the dimension $d$ increases; see Section~\ref{sec:comp_symb}.

The distinguishing feature of SMS is that they depend nonlinearly on their time-dependent parameters. This motivates the following definition.

\begin{defn}[SMS]
	We refer to $\hat u(\vc x,\btheta(t))$ as a shape-morphing solution of~\eqref{eq:gen_PDE} if it solves the PDE and is irreducibly nonlinear in $\btheta$. Here, \emph{irreducible nonlinearity} implies that for any two open sets $U, \tilde U\subset \mathbb R^n$, there exists no diffeomorphism $f:U\to \tilde U, \btheta\mapsto \tilde\btheta=f(\btheta)$ such that $\hat u(\vc x,f^{-1}(\tilde\btheta))$ is linear in $\tilde \btheta\in\tilde U$.
\end{defn}

\begin{rem}
The irreducible nonlinearity, instead of nonlinearity alone, is needed here to avoid idiosyncratic cases. For example, $\hat u(x,\theta(t))=x+\sin (t)$ is nonlinear in $\theta(t) = t$, but linear in $\tilde\theta(t) = \sin(t)$. Therefore, this function is not  irreducibly nonlinear in $\theta(t)=t$. In contrast, the function $\tanh(x+\sin(t))$ is irreducibly nonlinear in both $\theta(t)=t$ and $\tilde\theta(t)=\sin(t)$.
\end{rem}

We denote the set of all admissible parameter values by $\Omega$ so that $\btheta(t) \in\Omega\subseteq\mathbb R^n$.
If $\hat u(\vc x,\btheta(t))$ approximates the solution of the PDE (instead of being an exact solution), we refer to it as an \emph{approximate SMS}. However, since finding exact solutions of nonlinear PDEs are notoriously difficult, almost all solutions discussed in this review are approximate SMS, which we simply refer to as SMS.

\subsection{Examples of SMS}\label{sec:sms_ex}
In this section, we give a few concrete examples of shape-morphing solutions. 
\begin{enumerate}[wide,labelwidth=!,labelindent=0pt]
	\item \emph{Shape-morphing Gaussians:} These consist of a sum of Gaussian functions with time-dependent amplitudes $\alpha_i(t)\in\mathbb R$, centers or means $\pmb\mu_i(t)\in\mathbb R^d$, and lengthscales $\sigma_i(t)\in\mathbb R$,
	\begin{equation}
	\hat u(\vc x,\btheta(t))=\sum_{i=1}^r \alpha_i(t)\exp\left[ -\frac{|\vc x-\pmb\mu_i(t)|^2}{\sigma_i^2(t)}\right],
	\label{eq:sms_gauss_gen}
	\end{equation}
	where $\bbeta_i = \{\pmb\mu_i,\sigma_i\}$ are the shape parameters and $\alpha_i$ are the corresponding amplitudes. This results in 
	$n=(d+2)r$ set of parameters $\btheta = \{\alpha_i,\pmb\mu_i,\sigma_i\}_{i=1}^r$. The SMS~\eqref{eq:sms_gauss_gen} is similar to the kernels used in smoothed particle hydrodynamics~\cite{Monaghan1977,Lucy1977}.

	\item \emph{Shape-morphing (or evolutional) neural networks:} These are neural networks whose weights and biases are time-dependent. For instance, an evolutional shallow neural network can be written as 
	\begin{equation}
	\hat u(\vc x,\btheta(t))=\sum_{i=1}^r \alpha_i(t) \sigma(\vc w_i(t)\cdot \vc x+b_i(t)),
	\label{eq:sms_nn}
	\end{equation}
	where the shape parameters are $\bbeta_i=\{\vc w_i,b_i\}$, and $\btheta= \{\alpha_i,\vc w_i,b_i\}_{i=1}^r$. The sigmoidal function $\sigma:\mathbb R\to\mathbb R$ is the activation function of the network with internal weights $\vc w_i(t)$, biases $b_i(t)$, and output weights $\alpha_i(t)$. This construction has also been extended to deep neural networks~\cite{Zaki2021}. As we review in Section~\ref{sec:sms_eq}, the evolution of these parameters are prescribed by a set of ordinary differential equations (ODEs), and therefore no numerical optimization is needed to train the network. Consequently, evolutional neural networks do not suffer from generalization and training errors. We further discuss the shape-morphing neural networks in Section~\ref{sec:EDNN} and show that they can be readily generalized to deep neural networks.
	
	\item \emph{Shape-morphing finite elements:} As our last example, we introduce the shape-morphing finite elements in one spatial dimension,
	\begin{equation}
	\hat u(x,\btheta(t))=\sum_{i=1}^r \alpha_i(t) \phi_i(x;x_{i-1}(t),x_i(t),x_{i+1}(t)),
	\label{eq:sms_fem}
	\end{equation}
	where $x_i(t)$ denote the time-dependent grid points (or nodes) and
	\begin{equation}
	\phi_i(x;\cdots) = \begin{cases}
	\frac{x-x_{i-1}(t)}{x_i(t)-x_{i-1}(t)}, \quad x_{i-1}(t)\leq x\leq x_i(t),\\
	    		\\
	\frac{x_{i+1}(t)-x}{x_{i+1}(t)-x_{i}(t)}, \quad x_{i}(t)\leq x\leq x_{i+1}(t),\\
	    		\\
	0, \quad \mbox{otherwise,}
	\end{cases}
	\label{eq:sm_fem}
	\end{equation}
	is the shape-morphing hat function centered at the grid point $x_i(t)$. Unlike classical finite elements, where the grid points are independent of time, the shape-morphing finite elements~\eqref{eq:sms_fem} allow for the nodes to change over time to adapt to the solution of the PDE. In this example, the parameters are given by $\btheta(t) = [\alpha_1(t),x_1(t), \cdots, \alpha_r(t),x_r(t)]^\top\in\mathbb R^{2r}$. The grid points $x_0$ and $x_{r+1}$, at the boundary, are static.
\end{enumerate}

In each of these examples, we need to specify the evolution of parameters $\btheta(t)$ as a function of time. We review the optimal evolution equation for the parameters in Section~\ref{sec:sms_eq} below.

\subsection{Evolution of SMS}\label{sec:sms_eq}
In order to derive an evolution equation for the SMS parameters $\btheta(t)$, we substitute the SMS $\hat u(\cdot,\btheta(t))$ in the PDE~\eqref{eq:gen_PDE} and define the residual, 
\begin{equation}\label{eq:sms_res}
R(\vc x,\btheta,\dot\btheta) := \hat u_t - F(\hat u),
\end{equation}
where $\dot\btheta$ denotes the time derivative of the parameters $\btheta(t)$. The dependence of the residual on $\dot\btheta$ stems from the fact that
\begin{equation}
\hat u_t = \sum_{i=1}^n \pard{\hat u}{\theta_i}\dot{\theta}_i.
\end{equation}

If the SMS $\hat u$ is an exact solution, the residual vanishes. Otherwise, we evolve the parameters $\btheta(t)$ so that the residual $R$ is minimized in a suitable sense. 
For fixed $(\btheta,\dot\btheta)$, we view $R(\vc x,\btheta,\dot\btheta)$ as a function of $\vc x$, $R(\cdot,\btheta,\dot\btheta):\mathbb R^d\to\mathbb R$.
Assuming that this function belongs to a Hilbert space $\H$, we determine the evolution of parameters by solving the minimization problem,
\begin{equation}
\min_{\dot\btheta\in\mathbb R^n} \|R(\cdot,\btheta,\dot\btheta)\|_\H,
\label{eq:min_res}
\end{equation}
where $\|\cdot\|_\H$ is the norm induced by the Hilbert inner product $\langle\cdot,\cdot\rangle_\H$.

\begin{rem}\label{rem:DF_vari}
	Before discussing the solution to this minimization problem, a few technical remarks are in order.
	\begin{enumerate}[wide,labelwidth=!,labelindent=0pt]
		\item Equation~\eqref{eq:min_res} is an instantaneous optimization problem in the following sense. For a given set of parameters $\btheta(t)$ at time $t$, Eq. \eqref{eq:min_res} determines the rate of change of the parameters such that the norm of the residual at time $t$ is minimized. In quantum mechanics, this type of instantaneous optimization is called the Dirac--Frenkel variational principle~\cite{Beck2000,Raab2000}. As mentioned in Section~\ref{sec:history}, the same variational principle undergirds DLRA~\cite{Koch2007} and OTD~\cite{otd} methods in the context of matrix differential equations.
		
		\item 
		Instead of the instantaneous optimization, one could seek to minimize a finite-time notion of error, e.g., $\int_0^T \|R(\cdot,\btheta,\dot\btheta)\|^2_\H\id t$. However, as shown in the Appendix of \cite{Anderson2022}, the resulting Euler--Lagrange equations with Dirichlet boundary conditions can lead to spurious instabilities, and thus inaccurate evolution of SMS.
		
		\item Assuming that the residual $R$ belongs to a Hilbert space is largely for convenience, as it simplifies the forthcoming analysis. A similar optimization problem can be formulated if $R$ belongs to a Banach space. However, solving the corresponding minimization problem is less straightforward.
	\end{enumerate}
\end{rem}

In order to ensure that problem~\eqref{eq:min_res} is well-posed, we need to make the following assumption. 
\begin{ass}\label{ass:linInd}
	We assume the SMS $\hat u(\vc x,\btheta)$ is continuously differentiable with respect to its parameters $\btheta$ and that the partial derivatives with respect to these parameters are linearly independent. 
	That is the set, 
	\begin{equation}
	\left\{\pard{\hat u}{\theta_1}, \pard{\hat u}{\theta_2},\cdots,\pard{\hat u}{\theta_n} \right\},
	\end{equation}
	is a set of linearly independent functions of $\vc x$.
\end{ass}

Let's examine this assumption on the working example~\eqref{eq:ADE_sms} where $\theta_1=\alpha$, $\theta_2=\mu$, and $\theta_3=\sigma$ are the parameters. We have
\begin{subequations}\label{eq:working_ex_pd}
\begin{equation}
\pard{\hat u}{\theta_1} = \exp \left[ - \left(\frac{x-\mu}{\sigma}\right)^2\right],
\end{equation}
\begin{equation}
\pard{\hat u}{\theta_2} = \frac{2\alpha}{\sigma}\frac{x-\mu}{\sigma}\exp \left[ - \left(\frac{x-\mu}{\sigma}\right)^2\right],
\end{equation}
\begin{equation}
\pard{\hat u}{\theta_3} = \frac{2\alpha}{\sigma}\left(\frac{x-\mu}{\sigma}\right)^2\exp \left[ - \left(\frac{x-\mu}{\sigma}\right)^2\right].
\end{equation}
\end{subequations}
For any $\alpha\neq 0$, these three functions are linearly independent, which can be verified by computing their Wronskian that is equal to $(8\alpha^2/\sigma^5) \exp[-3(x-\mu)^2/\sigma^2]$.

Returning to the general case, the following result states that the minimization problem~\eqref{eq:min_res} has a unique solution and that the solution is known explicitly as an ODE for the parameters $\btheta$.
\begin{thm}
Under Assumption~\ref{ass:linInd}, the minimization problem~\eqref{eq:min_res} has a unique solution which satisfies, 
\begin{equation}\label{eq:sms_eq}
M(\btheta) \dot{\btheta} = \vc f(\btheta),
\end{equation}
where the metric tensor $M:\Omega\to\mathbb R^{n\times n}$ and the vector field $\vc f:\Omega\to\mathbb R^n$ are defined by
\begin{subequations}\label{eq:inn_prods}
\begin{equation}\label{eq:M}
M_{ij} := \left\langle \pard{\hat u}{\theta_i},\pard{\hat u}{\theta_j}\right\rangle_{\H}, 
\end{equation}
\begin{equation}\label{eq:rhs_vf}
f_{i} := \left\langle \pard{\hat u}{\theta_i},F(\hat u)\right\rangle_{\H},
\end{equation}
\end{subequations}
for $i,j\in\{1,2,\cdots,n\}$.
\end{thm}

The proof of this theorem can be found in \cite[Theorem 1]{Anderson2022}. We refer to \eqref{eq:sms_eq} as the \emph{SMS equation}.

\begin{rem}
Assumption~\ref{ass:linInd} ensures that the matrix $M(\btheta)$ is symmetric positive-definite, and hence invertible.
Therefore, the optimal evolution of parameters can be determined by solving the ODE,
\begin{equation}
\dot{\btheta} = M^{-1}(\btheta) \vc f(\btheta).
\end{equation}

If Assumption \ref{ass:linInd} is violated, the metric tensor might not be invertible. In that case, the inverse can be replaced by the Moore--Penrose pseudo inverse $M^+(\btheta)$, and the evolution of parameters is determined by the minimum-norm solution $\dot{\btheta} = M^+(\btheta)\vc f(\btheta)$. Since this approximation may introduce additional errors, it is best to choose the shape of the SMS judiciously, so that Assumption~\ref{ass:linInd} is satisfied.
\end{rem}

Returning to the working example, using the partial derivatives~\eqref{eq:working_ex_pd} and $\mathcal H=L^2(\mathbb R)$ as the Hilbert space, the metric tensor is given by 
\begin{equation}
M(\btheta) =\sqrt{\frac{\pi}{2}} \begin{bmatrix}
\sigma & 0 & \alpha/2\\
0 & \alpha^2/\sigma & 0\\
\alpha/2 & 0 & 3\alpha^2/4\sigma
\end{bmatrix},
\end{equation}
whose entries are computed by evaluating the integrals
\begin{equation}\label{eq:we_integrals}
M_{ij} = \int_{-\infty}^\infty \pard{\hat u}{\theta_i}(x,\btheta)\pard{\hat u}{\theta_j}(x,\btheta)\id x.
\end{equation}

Evaluating the right-hand side of the PDE~\eqref{eq:ADE} at the SMS~\eqref{eq:ADE_sms} yields,
\begin{align}
F(\hat u)  & = -c \hat u_x+\nu \hat u_{xx}	\nonumber \\
  &= \frac{2}{\sigma}  \left[ c \frac{x - \mu}{\sigma} + \frac{2\nu}{\sigma} \left(\frac{x-\mu}{\sigma}\right)^2 -  \frac{\nu}{\sigma}\right] \hat u.
\end{align}
Consequently, the right-hand side vector field in the SMS equation is 
\begin{equation}
\vc f(\btheta) =\sqrt{\frac{\pi}{2}} \begin{bmatrix}
-\alpha \nu/\sigma\\
\alpha^2 c/\sigma\\
\alpha^2\nu/2\sigma^2
\end{bmatrix}.
\end{equation}

Substituting these expressions in the SMS equation, and recalling that $\btheta = [\alpha,\mu,\sigma]^\top$, 
we obtain
\begin{equation}\label{eq:qdot_we}
	\dot \alpha = -\frac{2\alpha\nu}{\sigma^2},\quad \dot{\mu} = c,\quad \dot{\sigma}=\frac{2\nu}{\sigma},
\end{equation}
which admit the exact solution, 
\begin{equation}
	\alpha(t) = (4\pi\nu t)^{-\frac12}, \quad \mu(t) = ct, \quad \sigma(t) = (4\nu t)^{\frac12}.
\end{equation}
Substituting these expressions in the SMS~\eqref{eq:ADE_sms}, we observe that it coincides with the fundamental solution~\eqref{eq:ADE_FS} of the advection-diffusion equation.
In summary, the SMS equations dictate the correct parameter evolution for the SMS so that it coincides with the fundamental solution of the advection-diffusion equation.

\subsection{Enforcing conserved quantities}\label{sec:sms_eq_const}
Certain PDEs are derived from conservation laws such as conservation of mass, momentum, energy, etc~\cite{LeVeque_2002}.
Others admit secondary conserved quantities that are discovered only after the PDE is formulated~\cite{Lax1968}. 
In either case, an approximate solution should ideally also conserve these quantities. 
This is specially important for reduced-order models~\cite{Calberg2015} which do not approximate the solution accurately but need to 
preserve overall physical properties of the system. 

In the framework of shape-morphing solutions, it is straightforward to ensure that the ROM (or even a high-fidelity solution) respect the conserved quantities of the PDE. 
This is accomplished by adding the conserved quantities as constraints to the minimization problem~\eqref{eq:min_res}. 

\begin{defn}
	We refer to $I:\H\to\mathbb R$ as a conserved quantity (or a first integral) of the PDE~\eqref{eq:gen_PDE}, if $I(u(\cdot,t)) = I(u_0)$ for all times $t\geq 0$ and (almost) all initial conditions $u_0$.
\end{defn}

Evaluating a quantity $I$ along an SMS $\hat u(\cdot,\btheta(t))$, we can view the result as a function of the SMS parameters $\btheta(t)$.
As such, for notational simplicity, we write $I(\btheta(t))$ instead of  $I(\hat u(\cdot,\btheta(t)))$.

Our working example~\eqref{eq:ADE} does not admit a conserved quantity. But in the inviscid case $\nu = 0$, it is easy to verify that 
$I = \int_{-\infty}^{\infty} |u(x,t)|^2\id x$ is conserved. 
Evaluating this conserved quantity at the SMS~\eqref{eq:ADE_sms}, we obtain
\begin{equation}\label{eq:we_I}
I(\btheta) = \int |\hat u(x,\btheta)|^2\id x = \alpha^2\sigma \sqrt{\frac{\pi}{2}}.
\end{equation}

In general, consider a PDE that admits $m$ distinct conserved quantities $I_1,I_2,\cdots,I_m$.
To derive an evolution equation for the SMS parameters $\btheta(t)$, we add the conserved quantities $I_i$ as constraints to the 
minimization problem~\eqref{eq:min_res}. More precisely, we seek to solve the constrained Dirac--Frenkel variational principle,
\begin{subequations}\label{eq:min_res_const}
\begin{equation}
\min_{\dot\btheta} \|R(\cdot,\btheta,\dot{\btheta})\|_\H,
\end{equation}
\begin{equation}\label{eq:const}
\mbox{s.t.}\quad I_i(\btheta(t)) =\mbox{const.}, \quad i=1,2,\cdots,m.
\end{equation}
\end{subequations}

The solution to this constrained optimization problem is given by
\begin{equation}\label{eq:sms_eq_const}
	M(\btheta) \dot{\btheta} = \vc f(\btheta)- \sum_{k=1}^m \lambda_k \nabla I_k(\btheta),
\end{equation}
where $\lambda_k$ are Lagrange multipliers and $\nabla I_k$ denotes the gradient with respect to the parameters $\btheta$. We refer to~\eqref{eq:sms_eq_const} as the \emph{constrained SMS equation} and refer to \cite{Anderson2022} for its derivation.
The Lagrange multipliers $\pmb\lambda =[\lambda_1,\cdots,\lambda_m]^\top$ are obtained by solving the linear system, 
\begin{equation}\label{eq:lin_sys}
	C(\btheta) \pmb\lambda = \vc b(\btheta),
\end{equation}
where the constraint matrix $C(\btheta)\in\mathbb R^{m\times m}$ and the right-hand side vector $\vc b(\btheta)\in\mathbb R^m$ are given by 
\begin{equation}
	C_{ij} = \langle\nabla I_i,M^{-1}\nabla I_j\rangle,\quad b_i = \langle \nabla I_i,M^{-1}\vc f\rangle,
	\label{eq:Lag_lin_sys}
\end{equation}
where $\langle\cdot,\cdot\rangle$ denotes the usual Euclidean inner product and $\vc f$ is the right-side vector field~\eqref{eq:rhs_vf}.
If the gradients of the first integrals, $\{\nabla I_k\}_{k=1}^m$, are linearly independent, then the constraint matrix $C$ is symmetric positive-definite and therefore invertible.
Usually only a handful of conserved quantities are needed to be enforced and therefore solving the linear system~\eqref{eq:lin_sys} is not computationally costly.

For the working example in the inviscid case ($\nu=0$) and its first integral~\eqref{eq:we_I}, we have
\begin{equation}
	\nabla I(\btheta) = \sqrt{\pi/2}\, [2\alpha\sigma, 0, \alpha^2]^\top.
\end{equation}
The corresponding constraint matrix is given by
\begin{equation}
C = (\pi/2) \alpha^2(4\sigma^2+\alpha^2), \quad b=0.
\end{equation}
Note that $C$ and $b$ are scalars in this example since we only considered one constraint ($m=1$). Since the right-hand side vanishes, the resulting Lagrange multiplier is $\lambda=0$.
As a result, in this special case, the constraint SMS equation~\eqref{eq:sms_eq_const} coincides with the SMS equation. This is a consequence of the fact that the SMS equation~\eqref{eq:qdot_we} with $\nu=0$ happens to automatically conserve the first integral~\eqref{eq:we_I} and therefore no further enforcement is required.

\begin{rem}
In the constrained SMS~\eqref{eq:min_res_const}, one should be cautious not to over-constrain the problem. For instance, if the number of constraints $m$ surpasses the number of SMS parameters $n$, the constrained optimization problem may be ill-posed; in fact, the constraint set may even be empty. In other words, the model complexity $n$ should be large enough to allow for the existence of an SMS which also respects the conserved quantities.
\end{rem}
\begin{figure*}
	\centering
	\includegraphics[width=.95\textwidth]{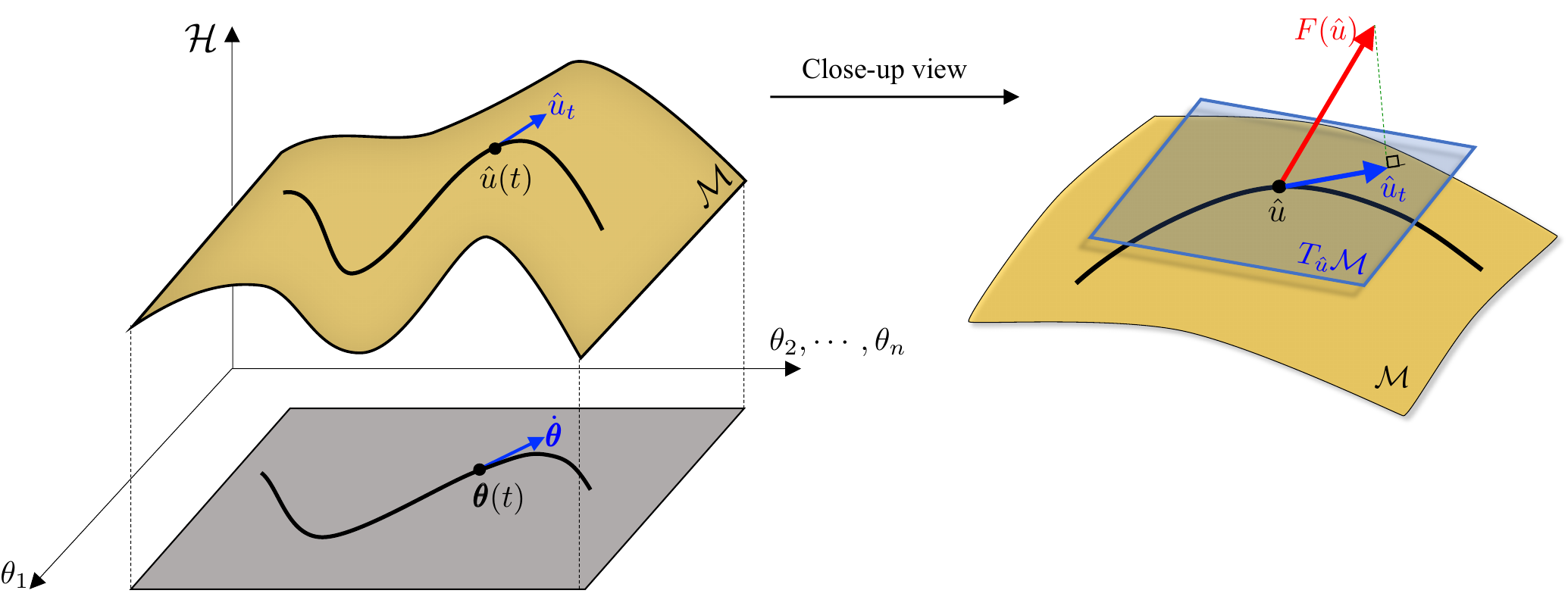}
	\caption{Illustrating the SMS manifold as a graph over its parameters. Here $\hat u(t)$ is shorthand for $\hat u(\cdot,\btheta(t))$.}
	\label{fig:schem_man}
\end{figure*}

\begin{rem}
Note that the constrained SMS equation~\eqref{eq:sms_eq_const} conserves the first integrals $I_k$ in continuous time. However, after time discretization, the conservation may be lost, unless a conservative discretization scheme is used. Alternatively, one can formulate and enforce the constraints at the discrete level as discussed by Schwerdtner et al.~\cite{Schwerdtner2024}.
\end{rem}

\section{Geometric interpretation of SMS}\label{sec:geom}
As discussed in~\cite{Anderson2022}, shape-morphing solutions have a nice geometric interpretation: SMS can be viewed as an $n$-dimensional submanifold of the functions space $\H$.

Consider the parameter set $\Omega$. An SMS can be viewed as a mapping 
$\hat u:\mathcal D\times \Omega\to\mathbb R$. For each fixed parameter $\btheta\in\Omega$, the SMS $\hat u(\cdot,\btheta)$ is a function of $\vc x\in\mathcal D$, belonging to the function space $\H$. The set of all such functions defines the \emph{SMS manifold},
\begin{equation}
\M := \bigcup_{\btheta\in\Omega} \hat u(\cdot,\btheta).
\label{eq:sms_man}
\end{equation}
As depicted in figure~\ref{fig:schem_man}, the SMS manifold can be viewed as a graph over the parameter space $\Omega$.
An evolution of parameters $\btheta(t)$ defines a curve in $\Omega$, which in turn defines a curve on the SMS manifold, corresponding to the SMS evolution $\hat u(\cdot,\btheta(t))$. In the special case of Galerkin approximation~\eqref{eq:Galerkin}, where $\hat u$ depends linearly on the amplitudes $\alpha_i$, the manifold is flat and spans a linear subspace of the Hilbert space $\H$.

The SMS manifold $\mathcal M$ is somewhat unusual: although it is an $n$-dimensional manifold, it lives in the infinite-dimensional function space $\H$; each point on this manifold is a function of $\vc x$. Nonetheless, it shares many properties of finite-dimensional manifolds. For instance, its tangent space is defined by
\begin{equation}
T_{\hat u} \M = \spn \left\{\pard{\hat u}{\theta_1},\pard{\hat u}{\theta_2}, \cdots, \pard{\hat u}{\theta_n}\right\}.
\end{equation}
Under Assumption~\ref{ass:linInd}, this tangent space is full-rank, meaning it is $n$-dimensional at all points on the manifold.
As such, $\M$ is an immersed submanifold of $\H$.

\begin{rem}
	The inner product $\langle\cdot,\cdot\rangle_\H$ on the Hilbert space $\H$ induces a Riemannian metric $\g$ on the parameter space $\Omega$. Consider any two vectors $\bxi,\etab\in\mathbb R^n$ at a point $\btheta\in\Omega$ and define
	\begin{align}
	\g_{\btheta} (\bxi,\etab) &:= \left\langle \sum_i\pard{\hat u}{\theta_i}\xi_i,\sum_j\pard{\hat u}{\theta_j}\eta_j\right\rangle_\H \nonumber\\
	& =\sum_{i,j}\xi_i\eta_j\left\langle \pard{\hat u}{\theta_i},\pard{\hat u}{\theta_j}\right\rangle_\H \nonumber\\
	& = \langle \bxi, M(\btheta)\etab\rangle,
	\end{align}
	where $\langle\cdot,\cdot\rangle$ denotes the usual Euclidean inner product and the last identity follows from~\eqref{eq:M}. Therefore, $\g$ is the pullback metric induced by the Hilbert inner product. This metric is identified with the matrix $M(\btheta)$, hence the name \emph{metric tensor}.
\end{rem}

The SMS equation~\eqref{eq:sms_eq} can be viewed as the orthogonal projection of the PDE dynamics on the tangent space of the SMS manifold $\M$.
More precisely, view the right-hand side of the PDE $F(\hat u)$ as a vector in $\H$ based at the SMS $\hat u$; see figure~\ref{fig:schem_man}. 
This vector does not necessarily belong to the tangent space $T_{\hat u}\M$. The SMS equation evolves the parameters $\btheta(t)$ such that the rate of change of the SMS $\hat u(\cdot,\btheta(t))$ coincides with the orthogonal projection of $F(\hat u)$ onto the tangent space $T_{\hat u}\M$. 
\begin{figure*}
	\centering
	\includegraphics[width=\textwidth]{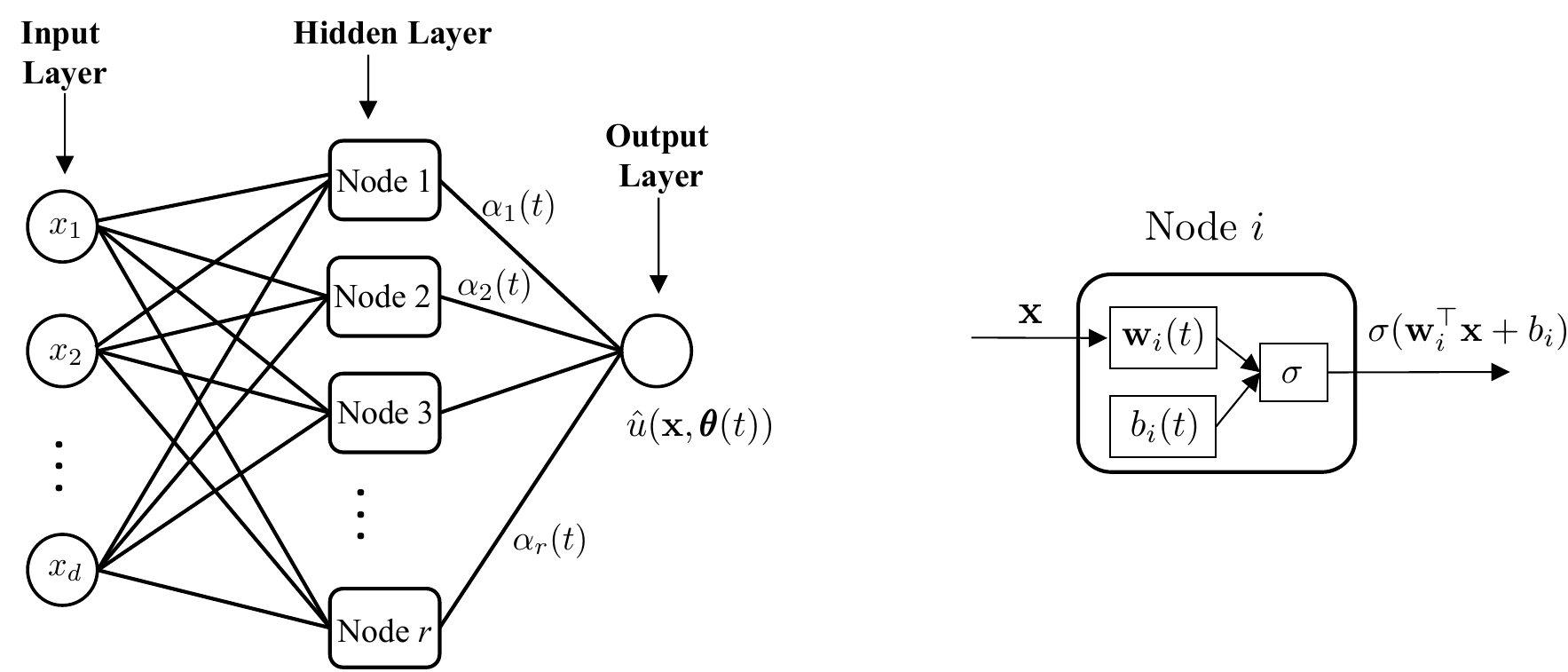}
	\caption{Architecture of a shallow evolutional neural network with $r$ nodes. The network parameters $\btheta(t)=\{\alpha_i(t),\vc w_i(t),b_i(t)\}_{i=1}^r$ are all time-dependent.}
	\label{fig:schem_net}
\end{figure*}

This geometric viewpoint immediately reveals alternative possibilities for determining the evolution equation for the parameters $\btheta(t)$. For example, instead of using the orthogonal projection which yields the SMS equation~\eqref{eq:sms_eq}, one can use an oblique projection onto the tangent space $T_{\hat u}\M$. 
Berman and Peherstorfer~\cite{Berman2023} used this geometric interpretation to devise their randomized sparse projection onto a submanifold of $\M$, which remedies overfitting and reduces computational cost. We will return to this work in Section~\ref{sec:comp_reg}.

\section{Neural networks as SMS}\label{sec:EDNN}
In general, the appeal of neural networks stems from their universal approximation property. For instance, Cybenko~\cite{Cybenko1989} showed that shallow neural networks with sigmoidal activation functions can approximate a continuous function with arbitrary accuracy. This result has been extended to deep neural networks and larger classes of functions such as Lebesgue spaces $L^p$ \cite{Park1991, Huang2006} and Sobolev spaces $W^{k,p}$ \cite{Hornik1991,Schmocker2025}. Therefore, it is natural to try to approximate the solution of a PDE using a neural network.

Although this viewpoint dates back to 1990s~\cite{Lagaris1998}, it has gained renewed momentum over the past decade with the introduction of physics-informed neural networks (PINNs)~\cite{Raissi2019,karniadakis2021} and operator learning methods~\cite{Stuart2020,Stuart2021,Lu2021}. These methods learn the static network parameters (i.e., weights and biases) through numerical optimization, e.g., stochastic gradient descent.

Du and Zaki~\cite{Zaki2021} introduced an interesting and alternative viewpoint. They envisioned the network parameters as time-dependent variables whose values change as the solution of the PDE evolves. They refer to such networks as evolutional neural networks, also known as neural Galerkin schemes~\cite{Bruna2024}.
For example, Eq.~\eqref{eq:sms_nn} describes an evolutional shallow neural network with time-dependent internal weights $\vc w_i(t)$, biases $b_i(t)$, and output weights $\alpha_i(t)$. Nonlinearities are encoded in the activation function $\sigma$. The network parameters therefore are $\btheta(t) = \{\alpha_i(t),\vc w_i(t),b_i(t)\}_{i=1}^r$ for a network with $r$ hidden nodes. Figure~\ref{fig:schem_net} shows the architecture of this neural network. 

This shallow architecture can easily be generalized to deep neural networks. Let $\vc N_\ell$ denote the output of the $\ell$-th layer of a deep neural network with $L$ hidden layers. The next layer is defined by the recursive relation
\begin{equation}\label{eq:DNN}
\vc N_{\ell+1}=\sigma(W_\ell(t)\vc N_\ell+\vc b_\ell(t)),
\end{equation}
where $\ell=0,1,\cdots, L$, the time-dependent matrix $W_\ell(t)$ determines the internal weights, and the vector $\vc b_\ell(t)$ denotes the biases. For the input layer, we have $\vc N_0 = \vc x$ and the output layer is given by $\hat u(\vc x,\btheta(t)) = W_{L+1}(t) \vc N_{L+1}+b_{L+1}(t)$. Collecting all the parameters of the network, we have $\btheta(t) = \{W_\ell(t), \vc b_\ell(t) \}_{\ell=0}^{L+1}$.

Clearly, evolutional neural networks are a special case of SMS. It is therefore not too surprising that Du and Zaki~\cite{Zaki2021}, in search of evolution equations for the parameters of the network, also arrived at the SMS equation~\eqref{eq:sms_eq}.

Compared to earlier neural network-based methods, evolutional neural networks have a number of advantages and an important drawback. Their main advantages are the following. First, since the evolution equations for the network parameters $\btheta(t)$ are explicitly known, evolutional neural networks do not require any data or training. Only the PDE being solved needs to be known. Second, because they do not require training, evolutional neural networks do not suffer from training or out-of-sample error.

On the other hand, evaluating an evolutional neural network requires solving a system of nonlinear ODEs, i.e., the SMS equation. 
As we discuss in Section~\ref{sec:comp}, this step can be time consuming. 
In contrast, although the training cost of PINNs and neural operators is high, their evaluation after training is very fast.

\section{Numerical examples}\label{sec:num_ex}
In this section, we present numerical examples where SMS are used to solve four different PDEs as listed in Table~\ref{tab:num_ex}: Nonlinear Schrodinger (NLS), Kuramoto–Sivashinsky (KS), Navier--Stokes (NS), and Fokker–Planck (FP) equations. For the first two examples (NLS and NS), we will present reduced-order (i.e., low fidelity) solutions of the PDEs, whereas the last two examples (KS and FP) present full-order (i.e., high fidelity) solutions.
\begin{table*}
	\centering
	\caption{Examples of PDEs solved with SMS: Nonlinear Schrodinger (NLS), Kuramoto–Sivashinsky (KS), Navier--Stokes (NS), Fokker–Planck (FP).}
	\begin{tabular}{|c|c|c|c|c|}
		\hline
		& \textbf{PDE} & \textbf{SMS} & \textbf{B.C.} & \textbf{Dimension}   \\ \hline
		\textbf{NLS} & $ \hat i u_t+u_{xx} + |u|^2u=0$ & Eq.~\eqref{eq:NLS_sms} & Decay & $d=1$  \\ \hline
		\textbf{NS} & $\omega_t+\vc v\cdot \nabla \omega =\nu\Delta \omega$ & Eq.~\eqref{eq:sms_gauss_gen} & Decay  & $d=2$  \\ \hline
		\textbf{KS} & $u_t+uu_x+u_{xx} + u_{xxxx}=0$ & Eq.~\eqref{eq:sms_KS} & Periodic & $d=1$  \\ \hline
		\textbf{FP} & $p_t+ \nabla\cdot(\vc v p)=\nu \Delta p$ & Eq.~\eqref{eq:sms_gauss_gen} & Decay & $d\geq 2$  \\ \hline
	\end{tabular}
	\label{tab:num_ex}
\end{table*}

Except for the KS equation which uses periodic boundary conditions, the other PDEs are defined on the infinite spatial domain with the assumption that the solutions (and their derivatives) decay to zero at infinity.

\subsection{Nonlinear Schr\"odinger equation}\label{sec:NLS}
The Nonlinear Schr\"odinger (NLS) equation is a perturbative approximation of dispersive waves appearing in water waves as well as optics~\cite{zakharov68}.
NLS is complex valued, $u(x,t)\in \mathbb C$, and its modulus $|u(x,t)|$ coincides with the wave envelope.
In particular, NLS has been used extensively to model and predict so-called rogue waves, i.e., waves of extreme amplitude~\cite{Onorato13,Farazmand2017a}.
Here, we show that SMS can be used to approximate rogue wave solutions of NLS in a reduced-order setting.

For the SMS we choose the Gaussian wave packet, 
\begin{equation}
\hat u(x,\btheta(t)) = \alpha(t) \exp\left[ -\frac{x^2}{\sigma^2(t)}+ \i\frac{x^2 v(t)}{\sigma(t)}+\i \varphi(t) \right],
\label{eq:NLS_sms}
\end{equation}
which was proposed by Desaix et al. \cite{Desaix91} to study the collapse of optical pulses. Here, $\i=\sqrt{-1}$, $\alpha(t)$ is the wave amplitude, and the shape parameters are $\{\sigma(t), v(t),\varphi(t)\}$ which denote the length scale, velocity, and phase of the wave, respectively.

There is no non-trivial exact solution of NLS that takes the form~\eqref{eq:NLS_sms}. However, certain solutions of NLS can be approximated by this SMS.
The corresponding SMS equation~\eqref{eq:sms_eq} for the parameters $\btheta = [\alpha(t),\sigma(t), v(t),\varphi(t)]^\top$ simplifies to
\begin{align}
\dot \alpha = -\frac{2\alpha v}{\sigma},\qquad & \dot\sigma = 4v\nonumber\\
\dot v = \frac{4}{\sigma^3}-\frac{\alpha^2}{\sqrt{2}\sigma},\qquad & \dot\varphi = \frac{5\alpha^2}{4\sqrt 2}-\frac{2}{\sigma^2}.
\label{eq:NLS_sms_eq}
\end{align}
This reduces the PDE to four coupled ODEs for the SMS parameters. Interestingly, the corresponding SMS~\eqref{eq:NLS_sms} automatically conserves mass and energy (Hamiltonian) of the system, although we did not directly impose these conservation laws in deriving equation~\eqref{eq:NLS_sms_eq}.

Figure~\ref{fig:num_ex}(a) shows the SMS starting from the initial conditions $\alpha(0) = 0.2$, $\sigma(0) = 20$, $v(0)=-0.05$, $\phi(0) = 0$. This mimics a rogue wave in the sense that the wave starts from a relatively small amplitude at the initial time $t=0$, grows to a relatively large amplitude (roughly twice the initial amplitude), and then decays back to a lower amplitude wave.
\begin{figure}
\centering
\includegraphics[width=0.55\textwidth]{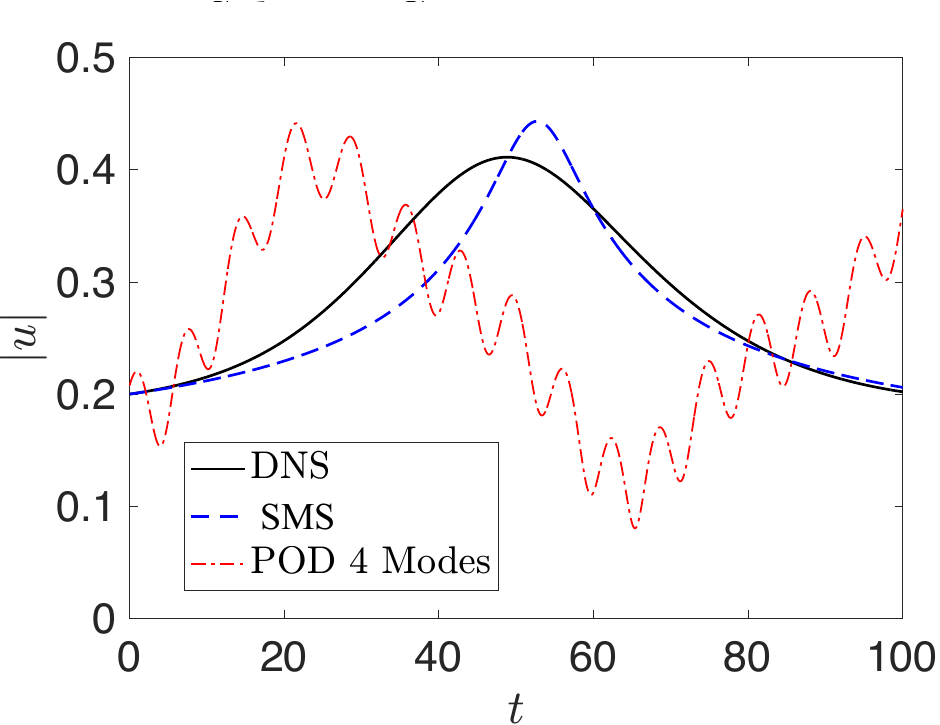}
\caption{Comparison for the NLS equation as originally reported in~\cite{Anderson2022}. The time series shows $|u(0,t)|$, the wave height at $x=0$. The SMS outperforms a comparable POD-based reduced-order model.}
\label{fig:NLS_compare}
\end{figure}

To compare this reduced-order SMS against a direct numerical solution (DNS) of NLS, we examine the amplitude of the wave as a function of time. This amplitude coincides with the modulus of the solution at $x=0$, i.e., $|u(0,t)|$. Figure~\ref{fig:NLS_compare} shows this comparison, highlighting that the SMS reproduces the amplitude and the timing of the maximum wave height with reasonable accuracy for a reduced-order solution. This figure also shows a ROM obtained from reducing NLS to four POD modes. We choose four POD modes to ensure fair comparison with the SMS which has four parameters. We observe that the SMS vastly outperforms the POD approximation.

It is noteworthy that SMS is an \emph{online} reduced-order modeling method. Whereas POD relies on training data to extract the POD modes in an offline phase, SMS does not require such training data. This offline training step is also required by many other ROMs such as Dynamic Mode Decomposition (DMD)~\cite{schmid10}. 
This online advantage is gained by choosing an appropriate shape for the SMS; Eq.~\eqref{eq:NLS_sms} in this example. As we discuss in Section~\ref{sec:concl}, a unified approach for choosing the optimal shape of SMS for general PDEs remains an open problem.

For further discussion on the application of SMS to dispersive waves, we refer to Anderson and Farazmand~\cite{Anderson2022b} who investigated SMS reduced-order solutions of NLS and modified NLS (a more accurate model of water waves). Among other things, they showed that the SMS solution remain accurate even when the wave groups propagate in space. Their theoretical results also revealed intricate connections between the SMS equation and the Lagrangian formulation of NLS.
Finally, we refer to Zhang et al.~\cite{Zhang2024} who applied SMS to the Gardner equation, a higher-order approximation of dispersive waves, in two spatial dimensions.

\subsection{Navier--Stokes equation}\label{sec:NS}
The Navier--Stokes equation, modeling fluid flow, has already been solved using SMS. For example, Du and Zaki~\cite{Zaki2021} use evolutional deep neural networks to solve the incompressible NS equation in two spatial dimensions with periodic boundary conditions. More recently, Kim and Zaki~\cite{Zaki2025} also solved the compressible Couette flow in 2D with no-slip boundary conditions.

\begin{figure*}
	\centering
	\includegraphics[width=0.95\textwidth]{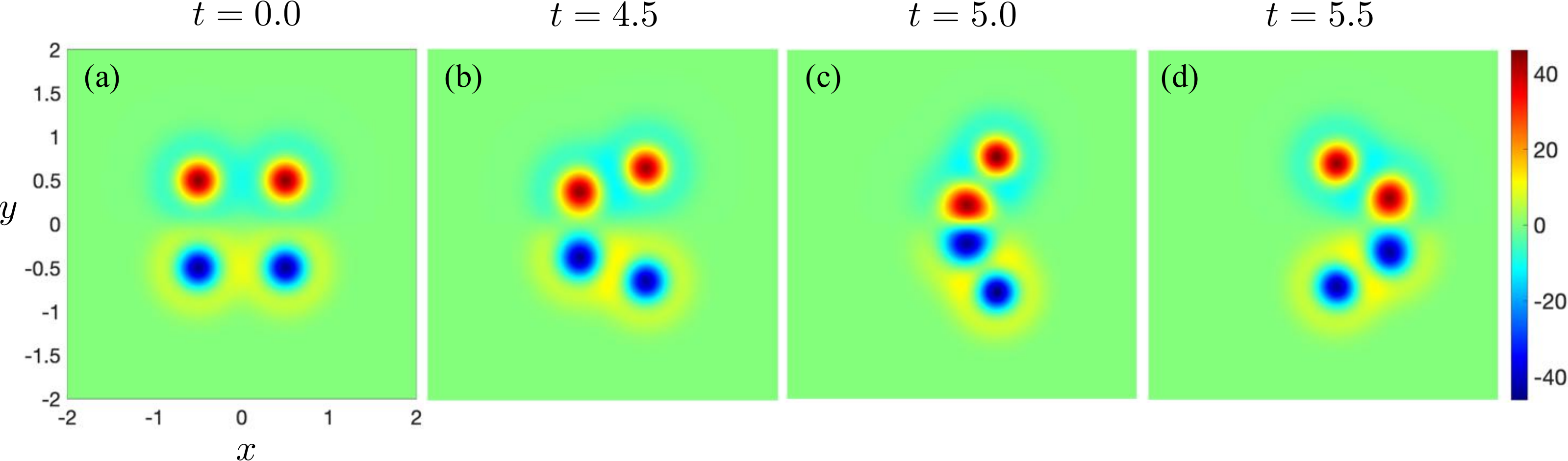}
	\caption{Leapfrogging dynamics of vortices reproduced by an SMS with $n=16$ parameters, as first reported in~\cite{Anderson2022}.}
	\label{fig:NS_leap}
\end{figure*}

Here, we showcase an application of SMS to the so-called leapfrogging phenomenon in vortex dynamics. Consider a 2D flow which is initialized as in figure~\ref{fig:NS_leap}(a): four vortices of the same strength but opposite signs arranged in a square configuration. It is well-known that this configuration exhibits a leapfrogging motion, whereby the trailing vortices speed up, deform, and pass in between the leading vortices. This motion repeats indefinitely.

We show that a simple SMS, consisting of four Gaussians reproduces this leapfrogging motion.
Consider the vorticity formulation of NS (see Table~\ref{tab:num_ex}) in the inviscid case $\nu=0$. In 2D, the stream function $\psi(\vc x,t)$ determines the 
velocity $\vc v$ and the vorticity $\omega$ fields so that $\vc v = [\partial_{x_2}\psi,-\partial_{x_1}\psi]^\top$ and $\omega = -\Delta \psi$.

For the SMS, we consider a stream function $\hat \psi(\vc x,\btheta(t))$ consisting of $r=4$ Gaussians as in \eqref{eq:sms_gauss_gen}.
The parameters of this SMS are the amplitudes $\alpha_i(t)\in\mathbb R$, lengthscales $\sigma_i(t)\in\mathbb R$ and the centers $\pmb\mu_i(t)\in\mathbb R^2$, resulting in $rK=16$ parameters $\btheta(t) = \{\alpha_i,\sigma_i,\pmb\mu_i\}_{i=1}^r$. These parameters are evolved according to the constrained SMS equation~\eqref{eq:sms_eq_const}. The first integrals enforced in these equations are the kinetic energy 
$I_1(\psi)= \int |\nabla \psi|^2\id \vc x$ and enstrophy $I_2(\psi)= \int |\Delta \psi|^2\id \vc x$.

Figure~\ref{fig:NS_leap} shows the resulting evolution of the SMS which qualitatively reproduces the leapfrogging motion. Anderson and Farazmand~\cite{Anderson2022} compare this evolution against point vortex dynamics
and report an excellent quantitative agreement as well.

To reproduce the leapfrogging vortex dynamics using SMS only a 16 dimensional ODE system (SMS equation) needs to be formed and solved.
It is worth pausing here and contemplate how this compares to conventional methods. To discretize the PDE, e.g. using finite differences or finite elements, one first needs to define an artificial boundary to contain the computations to a finite spatial domain. This domain needs to be large enough, with appropriate boundary conditions, to avoid unduly influencing the vortex dynamics~\cite{Hussain_2004}. Even then, to resolve the solution, the discretization will need to have far more degrees of freedom than the 16 parameters needed for SMS. Alternatively, imagine using POD modes to obtain a reduced-order solution. This approach would need to contend with two issues. First, a set of solutions of the vortex dynamics need to be obtained somehow from which the POD modes are extracted. Second, since the vortices advect to the right during their evolution, a set of globally defined static modes will struggle to accurately capture their dynamics. POD methods are generally known to flounder when applied to such advection-dominated flows~\cite{Reiss2018,Hilliard2024}.

\subsection{Kuramoto–Sivashinsky equation}\label{sec:KS}
The KS equation was originally derived to model reaction fronts. We consider this PDE on a spatial domain $x\in[-L/2,L/2]$ with $L=22$ and periodic boundary conditions. For the SMS, we use the evolutional shallow neural network~\eqref{eq:sms_nn} with a hyperbolic tangent activation function, $\sigma(\cdot) =\tanh(\cdot)$. An additional step needs to be taken to ensure that this SMS respects the periodic boundary conditions of the PDE. For example, we may apply a nonlinear change of coordinates to the spatial variable,
\begin{equation}
	\psi_i(x,t) = \sin \left( \frac{2\pi x}{L}+c_i(t)\right),
	\label{eq:sin_bc}
\end{equation}
and redefine the SMS to read, 
\begin{equation}\label{eq:sms_KS}
		\hat u(x,\btheta(t))=\sum_{i=1}^r \alpha_i(t) \sigma(w_i(t) \psi_i(x,t)+b_i(t)).
\end{equation}

Note that this SMS satisfies the periodic boundary conditions by construction. However, it comes at the expense of introducing the additional shift variables $c_i(t)$; without this shift, the network would fail to accurately approximate the KS solutions. Therefore, the SMS parameters are given by $\btheta = \{\alpha_i,w_i,b_i,c_i \}_{i=1}^r$, resulting in a total of $n = 4r$ parameters. We will further discuss the enforcement of boundary conditions in Section~\ref{sec:BC}.
\begin{figure}
	\centering
	\includegraphics[width=0.65\textwidth]{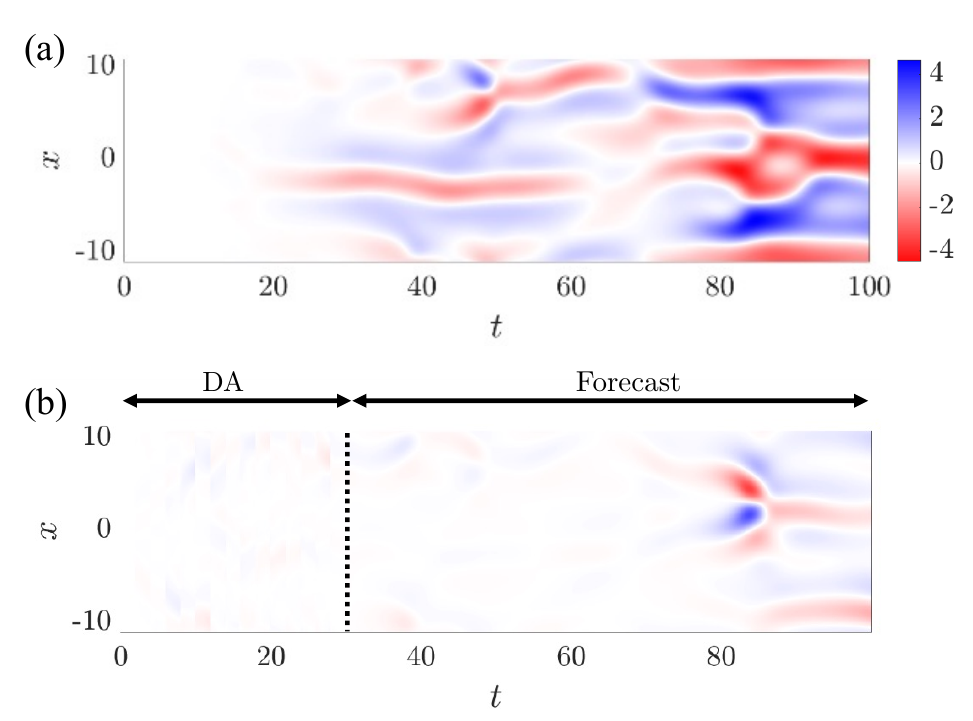}
	\caption{Error of the SMS for the KS equation, as first reported in~\cite{Hilliard2024b}. (a) Without data assimilation. (b) With data assimilation. In panel (b) data assimilation (DA) is only performed for the first 30 time units.}
	\label{fig:KS_error}
\end{figure}

Figure~\ref{fig:num_ex}(b) shows the corresponding SMS, with $r=10$ nodes, for the KS equation. The difference between this SMS and a DNS solution is shown in figure~\ref{fig:KS_error}(a). For the first 20 time units errors are small, but eventually the chaotic dynamics of KS amplifies the initially small errors. Data assimilation methods can be used to keep the errors small for a longer period of time; see figure~\ref{fig:KS_error}(b). We will return to data assimilation with SMS in Section~\ref{sec:concl}.

\subsection{Fokker--Planck equation}\label{sec:FP}
The Fokker--Planck equation describes the evolution of the probability density function (PDF) associated with a stochastic differential equation (SDE).
Consider the SDE, 
\begin{equation}
\id \vc X = \vc v(\vc X)\id t+ \sqrt{2\nu}\, \id \vc W,\quad \vc X(0) = \vc X_0,
\end{equation}
where $\vc v:\mathbb R^d\to\mathbb R^d$ is a sufficiently smooth vector field, $\vc W(t)\in \mathbb R^d$ is the standard Wiener process and $\nu>0$ is the diffusion coefficient. Because of the noise $\id\vc W$, the state of the system $\vc X(t)\in\mathbb R^d$ is also random and therefore needs to be understood probabilistically. The probability density $p(\vc x,t)$ of the state satisfied the FP equation in Table~\ref{tab:num_ex}. The initial condition $p(\vc x,0)=p_0(\vc x)$ is prescribed by the distribution of the initial state $\vc X_0$ which itself can be random. In the special case where the initial state $\vc X_0$ is known deterministically, the initial probability density is a Diract delta, $p_0(\vc x)=\delta(\vc x-\vc X_0)$.
\begin{figure}
	\centering
	\includegraphics[width=0.55\textwidth]{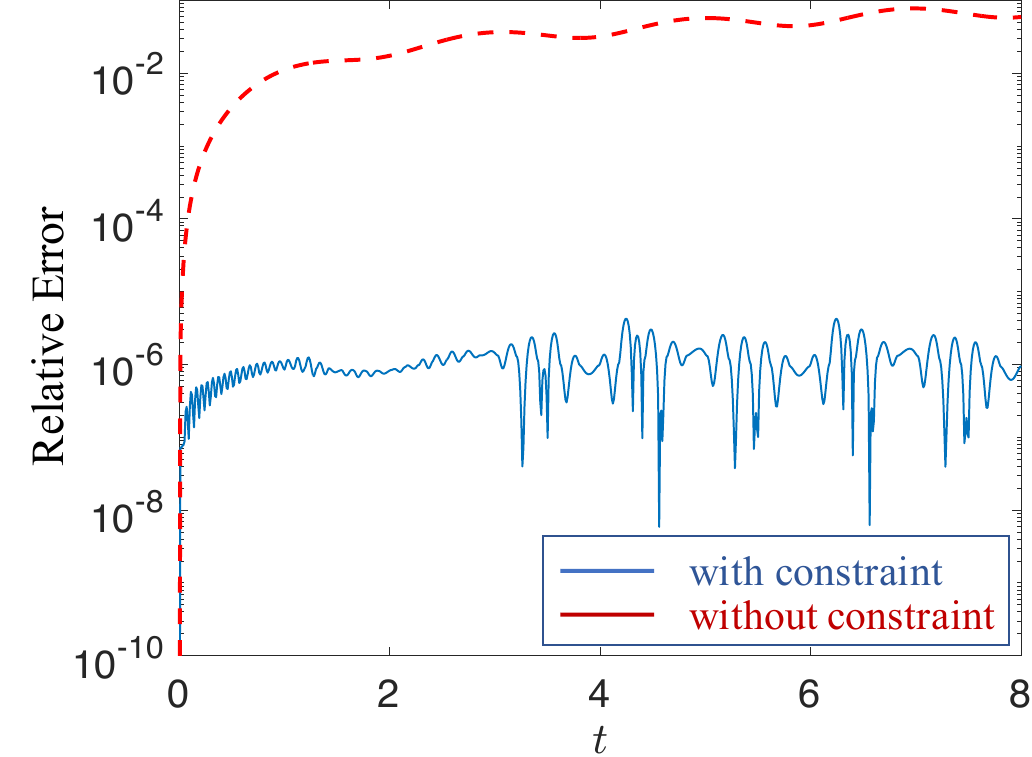}
	\caption{Enforcing the first integral~\eqref{eq:const_FP} in the SMS equation results in a more accurate approximation of the FP equation.}
	\label{fig:FP_err}
\end{figure}

As an example, we consider the stochastic Duffing oscillator, 
\begin{equation}
\ddot x = x-0.2\dot x-x^3 + \sqrt{2\nu} \dot W,
\label{eq:duffing}
\end{equation}
with $\nu =1/40$. Here, the state, $\vc X(t) = [x(t),\dot x(t)]^\top$, is two-dimensional ($d=2$). In the deterministic case ($\nu =0$), the damped Duffing oscillator~\eqref{eq:duffing} has two asymptotically stable fixed points $\vc X=(\pm 1,0)$ to which almost all trajectories converge. In the stochastic case, trajectories still tend to converge to a neighborhood of these equilibria such that the probability density $p(\vc x,t)$, in the asymptotic limit $t\to \infty$, is bimodal with its peaks near $(\pm 1,0)$. However, the transient dynamics can be quite complicated and depends on the initial density $p_0(\vc x)$.

Figure~\ref{fig:num_ex}(c) shows the evolution of an SMS corresponding to the FP equation of the stochastic Duffing oscillator. The SMS is a sum of $r=30$ axisymmetric Gaussians~\eqref{eq:sms_gauss_gen}. The initial condition $p_0(\vc x)$ is prescribed with $15$ Gaussians centered at $(1,1)$ and the remaining $15$ Gaussians at $(-1,-1)$. Figure~\ref{fig:num_ex}(c) shows this initial density and its progression through time according to the SMS equation. The Gaussians, which form the SMS, move and change shape until an asymptotic equilibrium is reached. By time $t=50$, the SMS has reached its asymptotic state.  Anderson and Farazmand~\cite{Anderson2024b} compared the SMS against large scale Monte Carlo simulations of the SDE~\eqref{eq:duffing}. They observed excellent agreement, not only for the asymptotic equilibrium density, but also for the transient dynamics of the density. In terms of computational time, it is notable that the Monte Carlo simulations take over two hours, whereas the SMS simulations take only about thirteen seconds.

In this example, the parameters of the SMS are evolved according to the constrained SMS equation~\eqref{eq:sms_eq_const}, where the first integral,
\begin{equation}
I(\btheta(t)) = \int_{\mathbb R^d} \hat p(\vc x,\btheta(t))\id \vc x =1,
\label{eq:const_FP}
\end{equation}
is enforced, ensuring that the SMS~\eqref{eq:sms_gauss_gen} is in fact a probability density function. 
Enforcing this constraint, in addition to ensuring that the results are physically meaningful, helps reduce the approximation error. 
For example, figure~\ref{fig:FP_err} shows the error for eight interacting particles in a harmonic trap (see \cite[Section 4.1]{Anderson2024a} for details). 
If constraint~\eqref{eq:const_FP} is enforced in the SMS dynamics, the relative error remains under $10^{-6}$. However, if the constraint is neglected, the error grows over time. 
This example highlights the importance of enforcing the conserved quantities of the PDE in the evolution equation of the corresponding SMS (see Section~\ref{sec:sms_eq_const}).

\begin{rem}[Fisher information metric]
	We conclude this section by pointing out an interesting feature of the SMS manifold~\eqref{eq:sms_man} in the case of the FP equation.
	In this case, the SMS manifold $\mathcal M$ is a statistical manifold, i.e., every point $\hat p(\cdot,\btheta)$ on it is a probability density parameterized by $\btheta$ \cite{Fisher1922}. 
	The intrinsic metric associated with a statistical manifold is 
	\begin{equation}
	\g_{ij}(\btheta) = \int_{\mathbb R^d} \frac{1}{\hat p(\vc x,\btheta)}\pard{\hat p}{\theta_i}(\vc x,\btheta)\pard{\hat p}{\theta_j}(\vc x,\btheta)\id\vc x,
	\label{eq:FIM}
	\end{equation}
	which is known as the \emph{Fisher information metric}.
	
	In the definition of the SMS metric tensor~\eqref{eq:M}, let the Hilbert space $\H$ to be a weighted $L^2$ space with the weight $[\hat p(\vc x,\btheta)]^{-1}$. Then the SMS metric tensor $M$ coincides with the Fisher information metric~\eqref{eq:FIM}. This reveals that, for the FP equation, the intrinsic Hilbert space is $\H = L^2_\mu (\mathbb R^d)$ where $\mu$ is an absolutely continuous measure with $\id \mu = [\hat p(\vc x,\btheta)]^{-1}\id \vc x$.
	
	It is notable that Bruna et al.~\cite{Bruna2024} estimated the inner products~\eqref{eq:inn_prods} using a Monte Carlo approximation. They proposed the above weighted inner product space $L^2_\mu$ heuristically as an importance sampling method to reduce the variance of the Monte Carlo estimates. Interestingly, in the special case of the FP equation, this heuristic is in fact rigorously justified as explained above.
\end{rem}

\section{Computational aspects}\label{sec:comp}
In this section, we review some computational bottlenecks of SMS as well as methods for addressing them.
The main contributors to the computational cost of SMS can be divided in two categories: (i) formation cost which refers to the computation of the metric tensor $M(\btheta)$ and the right-hand side vector field $\vc f(\btheta)$, and (ii) integration cost which refers to the numerical integration of the SMS equation in time.
Sections~\ref{sec:comp_symb} and \ref{sec:comp_coll} deal with the formation cost and Section~\ref{sec:comp_reg} addresses the integration cost.
In addition, we discuss the enforcement of boundary conditions of the PDE in Section~\ref{sec:BC}.

\subsection{Symbolic SMS}\label{sec:comp_symb}
Regarding the formation cost, note that evaluating the metric tensor $M$ and the right-hand side vector field $\vc f$ requires the computation of $\mathcal O(n^2)$ inner products~\eqref{eq:inn_prods} on the function space $\H$.
Each inner product requires the evaluation of an integral on the spatial domain $\mathcal D$; see, e.g., the working example in equation~\eqref{eq:we_integrals}.
Since these quantities depend on $\btheta(t)$, they need to be recomputed as the parameters evolve. The formation cost is not a major issue when SMS is used for reduced-order modeling since the number of parameters $n$ is relatively small. On the other hand, for high-fidelity solutions where the number of parameters $n$ is large, the formation cost becomes significant.

However, for a large class of SMS, the number of required inner product computations can be drastically reduced to $\mathcal O(K^2)$ where $K\ll n$. 
This is feasible if (i) the SMS takes the spectral form~\eqref{eq:sms_Galerkin}, and (ii) the inner products can be computed symbolically in closed form, e.g. using Mathematica or Maple.
All examples introduced in Section~\ref{sec:sms_ex} satisfy condition (i), but the evolutional \emph{deep} neural network~\eqref{eq:DNN} does not.
\begin{figure}
	\centering
	\includegraphics[width=0.55\textwidth]{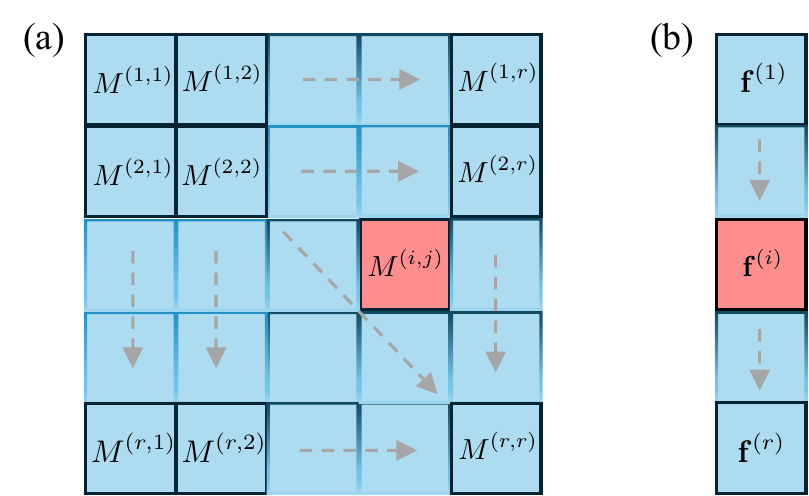}
	\caption{The partitioned matrix structure of the metric tensor (a) and the right-hand side vector field (b) for symbolic SMS.}
	\label{fig:symbSMS}
\end{figure}

If the above conditions are met, the SMS has $n=rK$ parameters, where $r$ is the number of terms in the summation \eqref{eq:sms_Galerkin} and
$K-1$ is the number of shape parameters in each term, $\bbeta_i =[\beta_{i1},\beta_{i2},\cdots,\beta_{i(K-1)}]^\top$. If the inner products in SMS are obtained via symbolic computing, we refer to the method as \emph{symbolic SMS}.

The resulting symbolic expressions allow us to only compute a fraction of the terms in the metric tensor and the right-hand side vector field. 
More specifically, as shown in figure~\ref{fig:symbSMS}(a), the metric tensor has a partitioned matrix structure with $K\times K$ blocks $M^{(i,j)}$.
Symmetry of the metric tensor dictates $M^{(i,j)}=M^{(j,i)}$. Furthermore, each block itself is symmetric, $[M^{(i,j)}]^\top = M^{(i,j)}$.
To evaluate all components of the metric tensor $M$, only a generic block $M^{(i,j)}$ needs to be computed symbolically.
The entries of this block consists of the inner products,
\begin{equation}
	\left\langle \pard{\hat u}{\alpha_i},\pard{\hat u}{\alpha_j}\right\rangle, \ \left\langle \pard{\hat u}{\alpha_i},\pard{\hat u}{\beta_{jk}}\right\rangle,
	\ \left\langle \pard{\hat u}{\beta_{ik}},\pard{\hat u}{\beta_{j\ell}}\right\rangle,
\end{equation}
where $k,\ell\in\{1,2,\cdots,K-1\}$. Once these symbolic expressions have been obtained, the entire metric tensor $M$ can be evaluated by inserting the numerical values of $\alpha_i$, $\alpha_j$, $\beta_{ik}$, etc.

A similar reduction in computational cost is applicable to the right-hand side vector field $\vc f:\mathbb R^n\to\mathbb R^n$. 
As shown in figure~\ref{fig:symbSMS}(b), this vector field has $r$ blocks of size $K\times 1$.
The components of the $i$-th block, $\vc f^{(i)}$, consist of the inner products,
\begin{equation}
	\left\langle \pard{\hat u}{\alpha_i},F(\hat u)\right\rangle_\H, \quad \left\langle \pard{\hat u}{\beta_{ik}},F(\hat u)\right\rangle_\H,
\end{equation}
where $k\in\{1,\cdots,K-1\}$. As long as this generic block is computed symbolically, the entire right-hand side vector field $\vc f$ can be evaluated by direct substitution of the numerical values of the corresponding parameters $\alpha_i$ and $\beta_{ik}$, etc.

In summary, evaluating the metric tensor $M$ and the vector field $\vc f$ require $K(K+3)/2$ symbolic computations. This leads to two major computational savings. First, the number of required symbolic computations is independent of the number of terms $r$ in the SMS. Second, the symbolic expressions need to be computed only once. As the time integration progresses, the updated numerical values of the parameters can be inserted in the previously computed symbolic expressions. 

The independence of the computational cost for symbolic SMS from the number of terms $r$ has significant consequences. As $r$ increases, so does the accuracy of the SMS approximation, but the formation cost remains constant. This allows symbolic SMS to be scalable for high-fidelity solutions. For example, the shape-morphing Gaussians~\eqref{eq:sms_gauss_gen}, which were used to solve the Navier--Stokes equation in Section~\ref{sec:NS} and the Fokker--Planck equation in Section \ref{sec:FP}, have three shape parameters ($K=4$) and therefore only require 14 symbolic computations, regardless of the number of Gaussians used. This property enabled \cite{Anderson2024b} to solve the Fokker--Planck equation accurately in higher spatial dimensions ($d=8$), which would have been impractical for conventional methods that rely on spatial discretization. We refer to \cite[Section 3.2]{Anderson2024a} for a more detailed discussion of symbolic SMS.

\subsection{Collocation SMS}\label{sec:comp_coll}
Obtaining symbolic expressions for the SMS equation is not always feasible. For example, Mathematica fails to return symbolic expressions for the metric tensor associated with the shallow neural network~\eqref{eq:sms_nn} with a hyperbolic tangent activation function.

In such cases, a collocation version of the SMS equations can be formulated which does not require any functional inner product computations. To describe this \emph{collocation SMS} method, consider the residual~\eqref{eq:sms_res}.
Instead of minimizing the norm $\|R(\cdot,\btheta,\dot\btheta)\|_\H$, we minimize the mean squared residual at prescribed collocation points. More specifically, let $\{\vc x_1,\vc x_2,\cdots,\vc x_N\}$ denote the collocation points. We seek to minimize the mean squared residual, 
\begin{equation}\label{eq:res_col}
\min_{\dot\btheta} \frac{1}{N} \sum_{i=1}^N |R(\vc x_i,\btheta,\dot\btheta)|^2. 
\end{equation}

As shown in \cite{Anderson2024a}, the minimum-norm optimizer is given by 
\begin{equation}\label{eq:col_sms_eq}
M^c(\btheta) \dot\btheta = \vc f^c(\btheta),
\end{equation}
where 
\begin{subequations}\label{eq:col_quant}
\begin{equation}\label{eq:col_M}
M^c_{ij}(\btheta) = \pard{\hat u}{\theta_j}(\vc x_i,\btheta), 
\end{equation}
\begin{equation}\label{eq:col_f}
f^c_i(\btheta)=F(\hat u(\cdot,\btheta))\vert_{\vc x=\vc x_i},
\end{equation}
\end{subequations}
for $i\in\{1,2,\cdots, N\}$ and $j\in\{1,2,\cdots, n\}$. The superscript here emphasizes that these quantities correspond to the collocation SMS. We refer to \eqref{eq:col_sms_eq} as the \emph{collocation SMS equation}.

Note that the collocation metric tensor $M^c$ is an $N\times n$ matrix. Viewing~\eqref{eq:col_sms_eq} as a linear equation for $\dot\btheta$, it may or may not have a solution. If it does have a solution, the residual $R$ vanishes at the collocation points along the corresponding SMS $\hat u(\cdot,\btheta(t))$. If it does not have a solution, the collocation SMS equation must be interpreted as a least squares problem which minimizes the residual~\eqref{eq:res_col}, and therefore, strictly speaking, is not a collocation point method. In either case, the evolution of parameters is given by the unique minimum-norm solution, $\dot\btheta = [M^c(\btheta)]^+\vc f^c(\btheta)$.

The advantage of the collocation SMS is that it can be formulated by the pointwise evaluations~\eqref{eq:col_quant}; it does not require any inner product (i.e., integral) evaluations. On the other hand, it may be less accurate since it only minimizes the residual at the collocation points instead of the residual over the entire spatial domain. Nonetheless, in practice, if the number of collocation points $N$ is large enough, the resulting SMS is quite accurate. For example, the KS equation in Section~\ref{sec:KS} was solved using collocation SMS with $N=128$.

\begin{rem}
	There is an intricate relationship between the SMS equation and the collocation SMS equation: 
	under certain conditions, the collocation SMS equation~\eqref{eq:col_sms_eq} is equivalent to the Monte Carlo approximation of the SMS equation~\eqref{eq:sms_eq}.

	For the collocation SMS, if $N>n$ and the columns of $M^c$ are linearly independent, then $[M^c]^+ = ([M^c]^\top M^c)^{-1}[M^c]^\top$.
	In this case, the collocation SMS equation is equivalent to the corresponding normal equation $[M^c(\btheta)]^\top M^c(\btheta) \dot\btheta = [M^c(\btheta)]^\top\vc f^c(\btheta)$.
	Next, note that 
	\begin{subequations}
		\begin{equation}
		\{[M^c(\btheta)]^\top M^c(\btheta)\}_{ij} = \sum_{k=1}^N \pard{\hat u}{\theta_i}(\vc x_k,\btheta)\pard{\hat u}{\theta_j}(\vc x_k,\btheta),
		\end{equation}
		\begin{equation}
		\{[M^c(\btheta)]^\top \vc f^c(\btheta)\}_{i} = \sum_{k=1}^N \pard{\hat u}{\theta_i}(\vc x_k,\btheta)F(\hat u(\cdot,\btheta))\vert_{\vc x=\vc x_k}.
	   \end{equation}
	\end{subequations}
These equations are exactly a Monte Carlo estimation of the inner products~\eqref{eq:inn_prods} with $\H = L^2_\mu(\mathcal D)$ (up to a scaling factor $1/N$) if the collocation points $\{\vc x_k\}$ are drawn at random according to the measure $\mu$.

In fact, Du and Zaki~\cite{Zaki2021} originally approximated the SMS equations using this Monte Carlo estimate. However, from the numerical standpoint, it is better to avoid the Monte Carlo estimate since the condition number of $[M^c(\btheta)]^\top M^c(\btheta)$ is the square of the condition number of $M^c(\btheta)$. Instead, to gain better numerical stability, one should compute the pseudo-inverse $[M^c(\btheta)]^+$ directly using singular value or QR decompositions~\cite{Golub1965,Golub2013}. 

\end{rem}

\subsection{Regularized SMS}\label{sec:comp_reg}
Now we turn to the computational challenges regarding the time integration of SMS equations. The (collocation)  SMS equation can be stiff as a set of nonlinear differential equations. As a result, explicit time integration schemes, e.g. Runge--Kutta, need to take exceedingly small time steps. This stiffness is related to the condition number of the metric tensor $M$.
When the number of SMS parameters $n$ is small, such as the reduced-order solutions discussed in Sections~\ref{sec:NLS}-\ref{sec:NS}, the metric tensor is often well conditioned and time integration is not an issue. However, as number of parameters increases, the metric tensor becomes more poorly conditioned and the resulting SMS equations exhibit hallmarks of stiffness.

An obvious resolution is to use an implicit integration scheme~\cite{Nick2025}. Although implicit schemes allow for larger time steps, they require solving a nonlinear system whose computational cost and convergence need to be taken into account. If one insists on explicit time integrators, a straightforward remedy was proposed in~\cite{Anderson2024a} by adding a regularization to the minimizing problem~\eqref{eq:min_res}. The regularized optimization problem reads, 
\begin{equation}\label{eq:res_reg}
	\min_{\dot\btheta} \|R(\cdot,\btheta,\dot\btheta)\|_\H + P(\dot\btheta),
\end{equation}
where $P:\mathbb R^n\to\mathbb R^+$ is the regularization (or penalty). For example, choosing the Tikhonov regularization, $P(\dot\btheta) = \gamma\|\dot\btheta\|^2$, the resulting \emph{regularized SMS equation} is given by 
\begin{equation}\label{eq:reg_SMS}
	\left( M(\btheta)+\gamma \mathbb I\right) \dot\btheta = \vc f(\btheta),
\end{equation}
where $\mathbb I$ denotes the identity matrix and $\gamma$ is the regularization parameter. 
This regularization idea can be similarly applied to the collocation SMS~\eqref{eq:res_col}.
We refer to~\cite{Anderson2024a} for a detailed derivation.

The regularized metric tensor $M+\gamma\mathbb I$ is much better conditioned than the metric tensor $M$. 
Consequently, explicit time integration, within reasonable computational time, becomes feasible. However, the regularized SMS 
introduces a free parameter, i.e., the regularization parameter $\gamma$, in the equations~\eqref{eq:reg_SMS}. As discussed \cite{Anderson2024a}, finding the optimal value of this regularization parameter is non-trivial.

Alternatives to direct regularization~\eqref{eq:res_reg} have also been explored. For example, Finzi et al. \cite{Finzi2023} use the randomized Nystr\"om preconditioning~\cite{Frangella2023} to achieve a better conditioned metric tensor $M$. This involves first obtaining a Nystr\"om approximation $M_{\mbox{nys}}$ of the metric tensor $M$ and a subsequent preconditioning $\mathbb P^{-1}M_{\mbox{nys}}$ of this approximation. The preconditioner $\mathbb P$ also contains a free parameter, a small regularization parameter whose optimal value needs to be specified.

Another alternative to regularization was proposed by Berman and Peherstorfer~\cite{Berman2023}. Motivated by the observation that the collocation metric tensor $M^c$ is low-rank, at each time step, they only update a subset of the parameters $\btheta$, leaving the rest unchanged. This is accomplished by defining a sparse parameter vector $\dot\btheta_s = S^\top\dot\btheta$, where $S$ is a sketching matrix, i.e., a subset of the columns of the identity matrix. More precisely, $S=[\vc e_{i_1}|\vc e_{i_2}|\cdots|\vc e_{i_k}]\in\mathbb R^{n\times k}$ where $\vc e_i$ are the standard basis on $\mathbb R^n$ and the indices $\{i_1,\cdots, i_k\}$ are a randomly selected subset of all possible indices $\{1,2,\cdots, n\}$ with $k\ll n$.

The collocation metric tensor is accordingly modified to $M^c(\btheta)S$. This procedure effectively modifies the collocation SMS equation~\eqref{eq:col_sms_eq} to $M^c(\btheta) SS^\top\dot\btheta = \vc f^c(\btheta)$. Note that this resembles a right preconditioning, although the orthogonal projection $SS^\top$ is not invertible. According to Berman and Peherstorfer~\cite{Berman2023}, this randomized projection approach remedies overfitting and reduces the computational cost. They test their method on three PDEs (Allen--Cahn, Burgers, and Vlasov equations) and report a few orders of magnitude reduction in the relative error.

\subsection{Boundary conditions}\label{sec:BC}
One of the major issues that remains to be addressed for SMS is the imposition of the boundary conditions of the PDE. Let's denote the boundary conditions by
\begin{equation}
	\B u(\vc x,t) = 0,\quad \forall \vc x\in\partial\mathcal D,
	\label{eq:BC_gen}
\end{equation}
where $\B$ is an appropriate boundary condition operator. For instance, for homogeneous Dirichlet boundary conditions, $\B$ is the restriction of $u$ to the boundary. For homogeneous Neumann boundary conditions, we have $\B u = \partial u/\partial n$, the directional derivative orthogonal to the boundary. Generalization to inhomogeneous boundary conditions is straightforward.

Unlike classical discretization methods, such as finite elements and finite differences, where boundary conditions can be easily enforced, a clear and rigorous methodology for imposing arbitrary boundary conditions for SMS remains illusive.

Nonetheless, some workarounds have been developed. For example, the idea of \emph{positional embeddings} was pursued by Du and Zaki~\cite{Zaki2021} as well as Kast and Hesthaven~\cite{Hesthaven2024}.  In this framework, one defines a map $\pmb\Psi : \mathcal D\subset\mathbb R^d\to \mathbb R^{d_e}$, where $d_e$ is the dimension of the embedding space. An intermediate SMS $U(\vc y,\btheta(t))$ is first defined on this embedding space where $\vc y\in\mathbb R^{d_e}$. Setting $\vc y =\pmb\Psi (\vc x)$, we obtain the final SMS, 
\begin{equation}
	\hat u(\vc x,\btheta(t)):=U(\pmb\Psi(\vc x),\btheta(t)) - U(\vc 0,\btheta(t)).
	\label{eq:sms_bc}
\end{equation}
For example, equation~\eqref{eq:sin_bc} is an example of such positional embedding for periodic boundary conditions.
In the context of evolutional neural networks (see Section \ref{sec:EDNN}), the positional embedding appears as an additional layer between the input layer and the first hidden layer of the network.

To ensure that the desired boundary conditions~\eqref{eq:BC_gen} are satisfied by the SMS~\eqref{eq:sms_bc}, the embedding $\pmb \Psi$ needs to satisfy certain conditions. In particular, consider the conditions,
\begin{subequations}\label{eq:BC_emb}
\begin{align}
	\mbox{Dirichlet:} & \quad \psi_i(\vc x)=  0, 	\label{eq:BC_emb_dir} \\
	\mbox{Neumann:}& \quad \nabla\psi_i(\vc x)\cdot \vc n(\vc x)= 0, \label{eq:BC_emb_neu}
	\end{align}
\end{subequations}
for all $\vc x\in\partial D$ and $i=1,2,\cdots, d_e$. Here, $\psi_i$ denote the components of $\pmb\Psi = [\psi_1,\cdots,\psi_{d_e}]$. 
If $\psi_i$ satisfy \eqref{eq:BC_emb_dir}, then the SMS \eqref{eq:sms_bc} satisfies the homogeneous Dirichlet boundary condition. 
Similarly, if $\psi_i$ satisfy \eqref{eq:BC_emb_neu}, then the SMS satisfies the homogeneous Neumann boundary condition. 

Obviously, the constant functions $\psi_i\equiv 0$ satisfy conditions~\eqref{eq:BC_emb}. However, they lead to a constant SMS~\eqref{eq:sms_bc}, oversimplifying the spatial dependence of the SMS. Therefore, the positional embeddings $\psi_i$ need to have enough variability inside $\mathcal D$ to obtain a non-trivial function $\hat u$. In order to achieve such embedding functions, in the context of Dirichlet boundary conditions, Du and Zaki~\cite{Zaki2021} used a function that depends on the distance of any point $\vc x\in\mathcal D$ to the boundary $\partial \mathcal D$. But their technique could not handle Neumann boundary conditions. 

Kast and Hesthaven~\cite{Hesthaven2024} pursued a more general framework which can tackle both Dirichlet and Neumann boundary conditions.
They proposed using the eigenfunctions of the Laplace operator,
\begin{align}
	\Delta \psi_i & = \lambda_i \psi_i, & \vc x\in\mathcal D,\nonumber\\
	\quad \B \psi_i & = 0, & \vc x\in\partial\mathcal D.
\end{align}
Other than limited cases, where these eigenfunctions are known analytically, they need to be obtained numerically, e.g., using a finite element method. This preprocessing step adds to the computational cost of the SMS and also impedes its applications in higher spatial dimensions. Furthermore, as Kast and Hesthaven~\cite{Hesthaven2024} point out, the embedding space dimension $d_e$ has to be large enough to ensure adequate spatial variability of the resulting SMS inside $\mathcal D$.

For rectangular domains, Hilliard and Farazmand~\cite{Hilliard2024b} introduced a domain doubling technique which allows for enforcement of 
Dirichlet, Neumann, and mixed Dirichlet-Neumann boundary conditions. We illustrate this method here for a one-dimensional domain $[0,L]$ and homogeneous Dirichlet boundary conditions. First we double the domain size by considering the enlarged domain $[-L,L]$, on which we define the periodic embedding, 
\begin{equation}
	\psi(x) = \sin\left( \frac{\pi x}{L}+c(t)\right).
\end{equation}
Note that this function is periodic over the enlarged domain $[-L,L]$. Next, we consider the SMS, 
\begin{equation}
	\hat u(x,\btheta(t)) = U(\psi(x),\hat\btheta(t))-U(\psi(-x),\hat\btheta(t)),
\end{equation}
where the continuous function $U(\cdot,\hat\btheta(t))$ is an intermediate SMS and $\btheta(t) = \{ \hat\btheta(t), c(t)\}$ are the parameters of the SMS.
It is easy to verify that $\hat u$ is $2L$-periodic and an odd function of $x$. These two properties ensure that $\hat u(0,\btheta(t))=\hat u(L,\btheta(t))=0$ for all times $t$. Therefore, the SMS automatically satisfies the Dirichlet boundary conditions on the original domain $[0,L]$. 

We refer to \cite[Appendix B]{Hilliard2024b} for generalization of this idea to higher dimensions as well as mixed Dirichlet-Neumann boundary conditions.
This domain doubling technique is limited to rectangular domains and suffers from the same curse of dimensionality as the positional embeddings.

We conclude this section by referring to Chen et al. \cite{Chen2023} who seek to enforce the boundary conditions by treating them as constraints. This is similar to 
enforcing the first integrals of the PDE in the constrained SMS equation (see Section~\ref{sec:sms_eq_const}).
However, this approach carries some subtle pitfalls that are discussed in \cite[Section D.1.2]{Hilliard2025_thesis}

\section{Conclusions and future directions}\label{sec:concl}
As we reviewed here, the theoretical and computational aspects of shape-morphing solutions have rapidly advanced over the past five years. In this section, we review the main open problems and computational challenges that remain to be addressed in the future.

\emph{Boundary conditions:} As reviewed in Section~\ref{sec:BC}, a few methods have been proposed for enforcing boundary conditions in SMS. 
However, the existing methods are far from being satisfactory or computationally optimal. For example, the positional embeddings require a user-specified embedding dimension $d_e$ with no clear theoretical guidelines for choosing this dimension. On the other hand, the domain doubling technique of \cite{Hilliard2024b} is not applicable to arbitrary geometries.
Furthermore, both methods are restricted to relatively low spatial dimensions ($d\leq 3$).
Therefore, developing a universal, unambiguous, and computationally efficient method for imposing boundary conditions remains an open problem.

\emph{Shape of SMS:} Perhaps the most nontrivial open problem is the choice of the SMS's shape for a given PDE. In the context of evolutional neural networks, the shape of the SMS is determined by the architecture and the activation functions of the network.
In Section~\ref{sec:NS}, we chose a sum of Gaussians to approximate vortex dynamics. This is a natural choice since vortex patches locally resemble a Gaussian~\cite{Trieling1997}. In contrast, the rogue wave SMS~\eqref{eq:NLS_sms}, which was motivated by previous studies on optical waves, is quite nontrivial. In general, the optimal choice for the shape of the SMS is not clear. Here, `optimal' can be interpreted as the shape that approximates the PDE, within a desired accuracy, with the smallest number of parameters $n$. This optimality is particularly important when SMS is used for reduced-order modeling, since ROMs need to have a minimal number of degrees of freedom. 

But even for high-fidelity solutions, the optimal shape is important. Although most studies do not clearly state this, solving PDEs in lower dimensions ($d\leq 3$) using SMS often remains computationally more costly than using conventional grid-based methods such as finite elements or finite differences. This is mainly due to the fact that a non-optimal choice of the SMS forces us to use too many parameters to satisfactorily approximate the PDE's solution, leading to unnecessarily large set of SMS equations.

Furthermore, a judiciously chosen SMS shape can reduce the stiffness of the SMS equations and the condition number of the metric tensor, hence reducing the computational cost of temporal integration (see Section~\ref{sec:comp_reg}). To date, no systematic method has been proposed for determining the optimal SMS shape for a given PDE.

\emph{Data assimilation with SMS:} In academic settings, one is usually content to solve a PDE from an arbitrary or user-defined initial condition $u_0$. However, if the PDE is being used to forecast the system that it models, the initial condition cannot be arbitrary; it must be inferred from the available observational data. 
This is the realm of data assimilation which is required for any numerical solver, including SMS. Hilliard and Farazmand \cite{Hilliard2024b} have developed a sequential data assimilation method which intermittently corrects the SMS parameters $\btheta(t)$ based on observational data. For example, figure~\ref{fig:KS_error}(b) shows the results of this data assimilation applied to the KS equation. The observations are pointwise measurements of the solution $u(x_i,t)$ at 10 equispaced sensors at the points $\{x_i\}_{i=1}^{10}$. This data is assimilated into the SMS solution for the first 30 time units. The SMS forecasts the system for the remaining 70 time units without data assimilation. As it is evident from the figure, the data assimilation significantly delays the growth of errors.

Although these results are promising, there is much more to be done on this front. For example, Kalman filtering, variational data assimilation, and nudging methods can in principle be used in conjunction with SMS. However, the nonlinear dependence of SMS on its parameters $\btheta$ complicates a straightforward application of these methods, hence requiring new theoretical and computational developments.

\emph{Weak solutions:} All shape-morphing solutions discussed here seek to approximate strong (sufficiently smooth) solutions of PDEs. However, if the SMS is a non-smooth function of space (e.g., the shape-morphing finite elements~\eqref{eq:sm_fem}) or if the PDE solution itself is non-smooth, a weak formulation of SMS needs to be used. Kim et al. \cite{kim2026} have begone work on weak formulation of SMS in the context of evolutional Kolmogorov--Arnold networks. 

Such a weak formulation can be obtained by requiring that $\langle R(\cdot,\btheta,\dot\btheta), v_k\rangle_\H =0$, i.e., the residual~\eqref{eq:sms_res} vanishes when evaluated against a set of test functions $v_k$. It remains to identify a set of test functions for a given SMS that ensure well-posedness of the inner products $\langle R(\cdot,\btheta,\dot\btheta), v_k\rangle_\H$. Furthermore, an optimally chosen set of test functions may lead to a sparse and well-conditioned metric tensor which will ameliorates the subsequent temporal integration of the weak SMS equations. These theoretical and computational questions remain to be addressed.

\emph{Error estimates:} In spite of the extensive body of work on SMS reviewed here, there are virtually no rigorous results that quantify their approximation error \cite{Gaby2024}. The universal approximation theorem is usually invoked to argue that as the number of parameters $n$ increases, the SMS will evermore closely approximate the true solution. Although this argument is correct, rigorous error estimates will specify the rate at which the error decays. 
Furthermore, error estimates will help avoid overfitting and may inform the choices we make about the shape of the SMS. This need is especially acute when SMS are used for model order reduction, where a minimal number of degrees of freedom is desirable. Hence, I conclude this review with an invitation to mathematical analysts to take this exciting and novel problem on.

\subsection*{Acknowledgments}
I am grateful to my former students William Anderson and Zack Hilliard who made significant contributions to the theory and computation of SMS. 
This work was supported by the National Science Foundation through the grant DMS-2208541.


\end{document}